\documentclass[11pt]{article}

\newcommand{\documentdate}{21 X 2021}

\usepackage{a4wide,latexsym,varioref,amsmath,amssymb,bbm,amsfonts,xcolor,graphicx}

\newcommand{\numsection}[1]{\section{#1}\setcounter{equation}{0}}
\newcommand{\appnumsection}[1]{\section*{#1}\setcounter{equation}{0}
  \renewcommand{\theequation}{A.\arabic{equation}}
  \renewcommand{\thetheorem}{A.\arabic{theorem}}
  \renewcommand{\thefigure}{A.\arabic{figure}}
  \renewcommand{\thesection}{A} }
\renewcommand{\theequation}{\arabic{section}.\arabic{equation}}
\renewcommand{\thefigure}{\arabic{section}.\arabic{figure}}

\newcommand{\calA}{{\cal A}}

\newcommand{\calE}{{\cal E}}
\newcommand{\calF}{{\cal F}}
\newcommand{\calG}{{\cal G}}
\newcommand{\calI}{{\cal I}}
\newcommand{\calM}{{\cal M}}
\newcommand{\calO}{{\cal O}}
\newcommand{\calS}{{\cal S}}

\newcommand{\calV}{{\cal V}}
\newcommand{\calZ}{{\cal Z}}
\newcommand{\beqn}[1]{\begin{equation}\label{#1}}
\newcommand{\eeqn}{\end{equation}}
\newcommand{\rr}{\overline{r}}
\newcommand{\ii}[1]{\{1, \ldots, #1 \}}
\newcommand{\iiz}[1]{\{0, \ldots, #1 \}}
\newcommand{\iibe}[2]{\{ #1, \ldots, #2 \}}
\newcommand{\flow}{f_{\rm low}}
\newcommand{\barphi}{\overline{\phi}}
\newcommand{\Df}{\Delta f}
\newcommand{\barDf}{\overline{\Delta f}}
\newcommand{\Dt}{\Delta t}
\newcommand{\barDt}{\overline{\Delta t}}
\newcommand{\barDT}{\Delta T}
\newcommand{\DF}{\Delta F}

\newcommand{\algn}{{\footnotesize {\sf TR$q$NE} }}
\newcommand{\expect}{\mathbb{E}}
\newcommand{\indic}{\mathbbm{1}_}
\newcommand{\no}[1]{#1^c}
\newcommand{\rlz}{(\omega)}
\newcommand{\AM}{\calM}
\newcommand{\tAM}{\widetilde{\calM}}
\newcommand{\indicI}{\indic{\calI}}
\newcommand{\indicM}{\indic{\AM_k}}
\newcommand{\indictM}{\indic{\widetilde{\AM}_k}}
\newcommand{\indicF}{\indic{\calF_k}}
\newcommand{\indictF}{\indic{\widetilde{\calF}_k}}
\newcommand{\indicE}{\indic{\calE_k}}
\newcommand{\indicEc}{\indic{\no{\calE_k}}}

\newcommand{\indicS}{\indic{\calS_k}}
\newcommand{\indicSc}{\indic{\no{\calS_k}}}
\newcommand{\prob}{\mathbb{P}{\rm r}}
\newcommand{\indicL}{\indic{\Lambda_k}}
\newcommand{\indicLc}{\indic{\no{\Lambda_k}}}
\newcommand{\indicbL}{\indic{\overline{\Lambda}_k}}

\DeclareMathOperator*{\argmax}{arg\,max}

\newcommand{\sigax}{\calA_{k-1}}
\newcommand{\pM}{p_*}

\newcommand{\sfrac}[2]{{\scriptstyle \frac{#1}{#2}}}
\newcommand{\half}{\sfrac{1}{2}}
\newcommand{\quarter}{\sfrac{1}{4}}
\newcommand{\tim}[1]{\;\; \mbox{#1} \;\;}
\newcommand{\mand}{\tim{ and }}
\newcommand{\sgn}{{\rm sgn}}
\newcommand{\req}[1]{(\ref{#1})}
\newcommand{\eqdef}{\stackrel{\rm def}{=}}
\newcommand{\ms}{\;\;\;\;}

\newcommand{\barf}{\overline{f}}
\newcommand{\bart}{\overline{t}}

\newcommand{\bigfrac}[2]{\frac{\displaystyle #1}{\displaystyle #2}}
\newcommand{\bigint}{\displaystyle \int}
\newcommand{\bigsum}{\displaystyle \sum}

\newtheorem{theorem}{Theorem}[section]
\newtheorem{lemma}[theorem]{Lemma}

\newtheorem{definition}{Definition}
\newcommand{\llem}[2]{\vspace{\baselineskip} 
\noindent\framebox[\textwidth]{\parbox{0.95\textwidth}{
\begin{lemma} \label{#1} \rm #2 \end{lemma} } } \vspace{\baselineskip} }
\newcommand{\lthm}[2]{\vspace{\baselineskip} 
\noindent\framebox[\textwidth]{\parbox{0.95\textwidth}{
\begin{theorem} \label{#1} \rm #2 \end{theorem} } } \vspace{\baselineskip} }
\newcommand{\bpr}{{\bf Proof.} \hspace{1.5mm}}
\newcommand{\epr}{\hfill $\Box$ \vspace*{1em}}
\newcommand{\proof}[1]{
\begin{list}{}{
\setlength{\topsep}{0.0pt}
\setlength{\partopsep}{0.0pt}
\setlength{\leftmargin}{0.025\textwidth}
\setlength{\rightmargin}{0.5\leftmargin}
\setlength{\labelwidth}{0.5\leftmargin}
\setlength{\labelsep}{0.25\leftmargin}}
\item \bpr #1 \epr \noindent
\end{list}}
\newcounter{algo}[section]
\renewcommand{\thealgo}{\thesection.\arabic{algo}}
\newcommand{\algo}[3]{\refstepcounter{algo}
\begin{center}\begin{figure}[htbp]
\framebox[\textwidth]{
\parbox{0.95\textwidth} {\vspace{\topsep}
{\bf Algorithm \thealgo : #2}\label{#1}\\
\vspace*{-\topsep} \mbox{ }\\
{#3} \vspace{\topsep} }}
\end{figure}\end{center}}
\newcommand{\al}[1]{{\footnotesize{\sf #1}}}
\newcommand{\ord}{\delta}
\newcommand{\ORD}{\Delta}
\newcommand{\sphi}{\widehat{\phi}}
\newcommand{\epsmin}{\epsilon_{\min}}
\renewcommand{\Re}{\mathbb{R}}
\newcommand{\neol}{\nonumber \\}

\newcommand{\comment}[1]{}

\title{
Trust-region algorithms: probabilistic complexity and intrinsic noise with applications to subsampling techniques
}

\author{
  S. Bellavia\thanks{Dipartimento di Ingegneria Industriale,
    Universit\`{a} degli Studi di Firenze, Italy. Member of the INdAM Research
    Group GNCS. Email: stefania.bellavia@unifi.it},
  G. Gurioli\thanks{Dipartimento di Matematica e Informatica ``Ulisse Dini'',
    Universit\`{a} degli Studi di Firenze, Italy.  {Member of the INdAM Research
    Group GNCS.} Email: gianmarco.gurioli@unifi.it},
  B. Morini\thanks{Dipartimento di Ingegneria Industriale,
    Universit\`{a} degli Studi di Firenze, Italy. Member of the INdAM Research
    Group GNCS. Email: benedetta.morini@unifi.it} \ and
  Ph. L. Toint\thanks{ Namur Center for Complex Systems (naXys),
    University of Namur, 61, rue de Bruxelles, B-5000 Namur, Belgium.
    Email: philippe.toint@unamur.be} 
}

\date{\documentdate}
\date{}

\begin{document}


\maketitle

\begin{abstract}
A trust-region algorithm is presented for finding approximate
minimizers of smooth unconstrained functions whose values and derivatives are
subject to random noise. It is shown that, under suitable probabilistic
assumptions, the new method finds (in expectation) an $\epsilon$-approximate minimizer of
arbitrary order $q\ge1$ in at most $\calO(\epsilon^{-(q+1)})$ 
inexact evaluations of the function and its derivatives,
providing the first such result for general optimality orders. 
The impact of intrinsic noise limiting the validity of the assumptions is also
discussed and it is shown that difficulties are unlikely to occur in
the first-order version of the algorithm for sufficiently large gradients.  Conversely,
should these assumptions fail for specific realizations, then ``degraded''
optimality guarantees are shown to hold when failure occurs.  These
conclusions are then discussed and illustrated in the context of subsampling
methods for finite-sum optimization.
\end{abstract}

 
\numsection{Introduction}

This paper is concerned with trust-region methods for solving the
unconstrained optimization problem
\beqn{problem}
\min_{x \in \Re^n} f(x),
\qquad f:\Re^n\rightarrow \Re,
\eeqn
where we assume that the \emph{values of the objective function $f$ and its
derivatives are computed subject to random noise}. Our objective is
twofold. Firstly, we introduce a version of the deterministic method proposed
in \cite{CartGoulToin20c} which is able to handle the random context and
provide, under reasonable probabilistic assumptions, a sharp evaluation
complexity bound (in expectation) for arbitrary optimality order. Secondly, we
investigate the effect of intrinsic noise (that is noise whose level cannot be
assumed to vanish) on a first-order version of our algorithm and prove
``degraded'' optimality, should this noise limit the validity of our
assumptions. The new results are then detailed and illustrated in the framework of
finite-sum minimization using subsampling.

Minimization algorithms using adaptive steplength and allowing for random
noise in the objective function or derivatives' evaluations have already
generated a significant literature
(e.g.\
\cite{BandScheVice14,ChenMeniSche18,BeraCaoSche19,PaquSche20,BlanCartMeniSche19,BellGuriMoriToin20,BellGuri21}).
We focus here on trust-region methods, which are methods in which a trial step
is computed by approximately minimizing a model of the objective function in a
``trust region'' where this model is deemed sufficiently accurate. The trial
step is then accepted or rejected depending on whether a sufficient  improvement in
objective function value predicted by the model is obtained or not, the radius
of trust-region being then reduced in the latter case.  We refer the reader to
\cite{ConnGoulToin00} for an in-depth coverage of this class of algorithms and
to \cite{Yuan15} for a more recent survey. Trust-region methods involving
stochastic errors in function/derivatives values were considered in particular in
\cite{BandScheVice14, CartSche17} and \cite{BlanCartMeniSche19,ChenMeniSche18}, 
the latter being the only methods (to the author's
knowledge) handling random perturbations in both the objective function and
its derivatives. The complexity analysis of the STORM (STochastic Optimization with Random Models)
algorithm described in
\cite{BlanCartMeniSche19,ChenMeniSche18} is based on supermartingales and makes probabilistic
assumptions on the accuracy of these evaluations which become tighter when the
trust-region radius becomes small.  It also hinges on the definition of a
monotonically decreasing ``merit function'' associated with the stochastic
process corresponding to the algorithm. The method proposed in this paper can
be viewed as an alternative in the same context, but differs from the STORM
approach in several aspects. 
The first is that the method discussed here uses a model whose degree is
chosen adaptively at each iteration, requiring the (noisy) evaluation of
higher derivatives only when necessary. The second is that its scope is not
limited to searching for first- and second-order approximate minimizers, but
is capable of computing them to arbitrary optimality order. The third is that
 the probabilistic accuracy
conditions on the derivatives' estimations no longer depends on the trust-region
radius, but rather on the predicted reduction in objective function
values, which may be less sensitive to problem conditioning. Finally,
its evaluation complexity analysis makes no use of
a merit function of the type used in \cite{BlanCartMeniSche19}.

In \cite{BellGuriMoriToin21b}, the impact of intrinsic random noise on the
evaluation complexity of a \emph{deterministic} ``noise-aware'' trust-region
algorithm for unconstrained nonlinear optimization was
investigated and constrasted with that of an inexact version where noise is
fully controllable. The current paper considers the question in the more general
probabilistic framework.

Even if the analysis presented below does not depend in any way on the
choice of the optimality order $q$, the authors are well aware that, while
requests for optimality of orders $q\in\{1,2\}$ lead to
practical, implementable algorithms, this may no longer be the case for
$q>2$. For high orders, the methods discussed in the paper therefore constitute an
``idealized'' setting (in which complicated subproblems can be approximately
solved without affecting the evaluation complexity) and thus indicate the
limits of achievable results.

The paper is organized as follows. After introducing the new stochastic trust-region algorithm
in Section~\ref{section:algorithm}, its evaluation complexity
analysis is presented in Section~\ref{section:analysis}.
Section~\ref{noise-s} is then devoted to an in-depth discussion of the impact
of noise on the first-order instantiation of the algorithm, with a particular
emphasis on the case where noise is generated by subsampling in finite-sum
minimization context. Conclusions and perspectives are
finally proposed in Section~\ref{section:conclusion}.
Because our contribution borrows ideas from \cite{BellGuriMoriToin20},
themselves being partly inspired by \cite{CartSche17}, repeating some material
from these sources is necessary to keep our argument understandable. We have
however done our best to limit this repetition as much as possible.

{\bf Basic notations.} Unless otherwise specified, $\|\cdot\|$ denotes the standard
Euclidean norm for vectors and matrices.  For a general symmetric tensor $S$
of order $p$, we define
\[
\|S\|_{[p]} \eqdef \max_{\|v\|=1}  | S [v]^p |
= \max_{\|v_1\|= \cdots= \|v_p\|=1} | S[v_1, \ldots, v_p] |
\]
the induced Euclidean norm. We also denote by $\nabla_x^j f(x)$ the $j$-th
order derivative tensor of $f$ evaluated at $x$ and note that such a tensor is
always symmetric for any $j\geq 2$. $\nabla_x^0 f(x)$ is a synonym for $f(x)$.
$\lceil \alpha \rceil$ denotes the
smallest integer not smaller than $\alpha$. Moreover, given a set ${\cal B}$, $|{\cal B}|$ denotes its cardinality, $\indic{\cal B}$
refers to its indicator function and $\no{{\cal{B}}}$ indicates its complement.
 All stochastic quantities live in a probability space denoted by  $(\Omega, \calA, \prob)$ with the probability measure $\prob$ and the $\sigma$-algebra $\calA$
 containing subsets of $\Omega$. We never explicitly define $\Omega$, but specify it through random variables. 
 $\prob[{\rm event}]$ finally denotes the probability of an
event and $\expect[X]$ the expectation of a random variable $X$. 

\numsection{A trust-region minimization method for problems with \\randomly 
  perturbed function values and derivatives}\label{section:algorithm}

We make the following assumptions on the optimization problem \req{problem}.
\vskip 5pt
\begin{description}
\item[AS.1] The function $f$ is $q$-times continuously
  differentiable in $\Re^n$, for some $q\ge1$. Moreover, its $j$-th
  order derivative tensor is  Lipschitz continuous for $j\in\ii{q}$ in the
  sense that, for each $j\in\ii{q}$, there exists a constant $\vartheta_{f,j}\geq 0$
  such that, for all $x,y\in \Re^n$,
  \beqn{f-holder}
  \|\nabla_x^j f(x) - \nabla_x^j f(y)\| \leq \vartheta_{f,j} \|x-y\|.
  \eeqn
\item[AS.2] $f$ is bounded below in $\Re^n$, that is there exists a constant
  $\flow$ such that $f(x)\ge \flow$ for all $x \in \Re^n$.
\end{description}
\noindent
\vskip 5pt \noindent

\noindent
AS.2 ensures that the minimization problem \req{problem} is well-posed. 
AS.1 is a standard assumption in evaluation complexity analysis\footnote{It
is well-known that requesting \req{f-holder} to hold for
all $x,y\in \Re^n$ is strong. The weakest form of AS.1 which we could
use in what follows is to require \req{f-holder} to hold for all $x = x_k$
(the iterates of the minimization algorithm we are about to describe) and all
$y= x_k+ \xi s_k$ (where $s_k$ is the associated step and $\xi$ is arbitrary
in [0,1]). However, ensuring this condition a priori, although maybe possible for
specific applications, is hard in general, especially for a non-monotone algorithm with a
random element.}. It is important because we
consider algorithms that are able to exploit all available derivatives of $f$
and, as in many minimization methods, our approach is based on
using the Taylor expansions 
(now of degree $j$ for $j\in\ii{q}$) given by
\beqn{taylor}
t_{f,j}(x,s) \eqdef f(x) + \sum_{\ell=1}^j \nabla_x^\ell f(x)[s]^\ell.
\eeqn
AS.1 then has the following crucial consequence.

\llem{tech-Taylor-theorem}{
Suppose that AS.1 holds. Then for all $x,s \in \Re^n$,
\beqn{tech-resf}
|f(x+s) - t_{f,j}(x,s)| \leq \bigfrac{\vartheta_{f,j}}{(j+1)!} \, \|s\|^{j+1}.
\eeqn
}
\proof{See \cite[Lemma~2.1]{CartGoulToin20b} with $\beta=1$.}
\vskip 5pt
\noindent

\noindent
At a given iterate $x_k$ of our algorithm, we will be interested in finding a
step $s\in \Re^n$ which makes the \emph{Taylor decrements} 
\beqn{DT-def}
\Dt_{f,j}(x_k,s) \eqdef f(x_k)- t_{f,j}(x_k,s) =
t_{f,j}(x_k,0)-t_{f,j}(x_k,s)
\eeqn
large (note that $\Dt_{f,j}(x,s)$ is independent of $f(x)$). When this is
possible, we anticipate from the approximating properties of the Taylor
expansion that some significant decrease is also possible in $f$.  Conversely,
if $\Dt_{f,j}(x,s)$ cannot be made large in a neighbourhood of $x$, we
must be close to an approximate minimizer.  More formally, we define, for some
$\theta\in (0,1]$ and some \emph{optimality radius} $\ord \in (0,\theta]$,
the measure 
\beqn{phi-def}
\phi_{f,j}^{\ord}(x) = \max_{\|d\|\leq \ord} \Dt_{f,j}(x,d),
\eeqn
that is the maximal decrease in $t_{f,j}(x,d)$ achievable in a ball of radius $\ord$
centered at $x$.  (The practical purpose of introducing $\theta$ is to avoid
unnecessary computations, as discussed below.) We then define $x$ to be a $q$-th order
$(\epsilon,\ord)$-approximate minimizer (for some accuracy requests
$\epsilon\in (0,1]^q$) if and only if  
\beqn{approx-min}
\phi_{f,j}^\ord(x) \leq \epsilon_j \frac{\ord^j}{j!}
\tim{ for } j \in \ii{q},
\eeqn
(a vector $d$ solving the optimization problem defining
$\phi_{f,j}^{\ord}(x)$ in \req{phi-def} is called an \emph{optimality displacement}) \cite{CartGoulToin17,CartGoulToin20c}.
In other words, a $q$-th order $(\epsilon,\ord)$-approximate
minimizer is a point from which no significant decrease of the Taylor
expansions of degree one to $q$ can be obtained in a ball of optimality radius
$\ord$. This notion is coherent with standard optimality measures for low
orders\footnote{It is easy to verify that, irrespective of $\ord$,
\req{approx-min} holds for $j=1$ if and only if $\|\nabla_x^1f(x)\|\leq\epsilon_1$
and that, if $\|\nabla_x^1f(x)\|=0$, $\lambda_{\min}[\nabla_x^2 f(x)] \geq -\epsilon_2$ if
and only if $\phi_{f,2}^\ord(x) \leq \half \epsilon_2\ord^2$.} and has the advantage of
being well-defined and continuous in $x$ for every order.
Note that $\phi_{f,j}^\ord(x)$  is  always non-negative.

This paper is concerned with the case where the values of the objective
function $f$ and of its derivatives $\nabla_x^j f$  are subject to random
noise and can only be computed inexactly (our assumptions on random noise will
be detailed below). Our notational convention will be to 
denote inexact quantities with an overbar, so $\barf(x,\xi)$ and
$\overline{\nabla_x^j f}(x,\xi)$ denote inexact values of $f(x)$ and $\nabla_x^j
f(x)$, where $\xi$ is a random variable causing inexactness. Thus \req{taylor}
and \req{DT-def} are unavailable, and we have to consider
\[
\bart_{f,j}(x_k,s,\xi) \eqdef \barf(x_k,\xi) + \sum_{\ell=1}^j \overline{\nabla_x^\ell f}(x_k,\xi)[s]^\ell
\]
and the associated decrement
\beqn{barDT-def-xi}
\barDt_{f,j}(x,s_k,\xi)
\eqdef \bart_{f,j}(x_k,0,\xi)-\bart_{f,j}(x_k,s_k,\xi)
= -\sum_{\ell=1}^j \overline{\nabla_x^\ell f}(x_k,\xi)[s]^\ell
\eeqn
instead. For simplicity, we will often omit to mention the dependence of inexact
values on the random variable $\xi$ in what follows, so \req{barDT-def-xi}
is rewritten as
\beqn{barDT-def}
\barDt_{f,j}(x,s_k)
\eqdef \bart_{f,j}(x_k,0)-\bart_{f,j}(x_k,s_k)
= -\sum_{\ell=1}^j \overline{\nabla_x^\ell f}(x_k)[s]^\ell.
\eeqn
This in turn would require that we measure optimality using
\beqn{phi-def-bar}
\barphi_{f,j}^{\ord}(x) \eqdef
\max_{ \|d\| \leq \ord} \barDt_{f,j}(x,d)
\eeqn
instead of \req{phi-def}. However, computing this exact global maximizer may be
costly, so we choose to replace the computation of \req{phi-def-bar} by an
approximation, that is with the computation of an optimality displacement
$d$ with $\|d\|\le \ord$ such that
$
\varsigma \barphi_{f,j}^{\ord}(x) \le \barDt_{f,j}(x,d)
$ 
for some constant $\varsigma \in (0,1]$. 
We  state the Trust-Region with Noisy Evaluations (\algn) algorithm  \vpageref{algoS} 
using all the ingredients we have described. The trust region radius at iteration $k$ is denoted by $r_k$
instead of the standard notation $\Delta_k$. 

\algo{algoS}{The \algn algorithm}
{
\begin{description}
\item[Step~0: Initialisation.]
  A criticality order $q$, a starting point $x_0$ and accuracy
  levels $\epsilon \in (0,1)^q$ are given. For a given constant $\eta
  \in (0,1)$, define
  \beqn{epsmin-omega-def}
  \epsilon_{\min}\eqdef \min_{j\in\{1,..,q\}}\epsilon_j
  \mand
  \nu \eqdef \min\big[\half \eta, \quarter(1-\eta)\big].
  \eeqn
  The constants $\theta \in [\epsilon_{\min}, 1]$, $\varsigma \in
  (0,1]$,  $\gamma>1$, $r_{\max}\ge 1$ and an initial trust-region radius
  $r_0\in (\epsilon_{\min}, r_{\max}]$ are also given. Set $k=0$.

\vskip 5pt
\item[Step~1: Derivatives estimation.]
  Set $\ord_k = \min[r_k,\theta]$. For $j = 1, \ldots, q $,
  \vspace*{-2mm}
  \begin{enumerate}
  \item Compute derivatives' estimates $\overline{\nabla_x^j f}(x_k)$ and
    find an optimality displacement $d_{k,j}$ with $\|d_{k,j}\|\le\ord_k$ such that
 \begin{equation}
 \label{globmax_inexact}
 \varsigma \barphi_{f,j}^{\ord_k}(x_k) \le \barDt_{f,j}(x_k,d_{k,j}).
  \vspace*{-2mm}
 \end{equation}
  \vspace*{-2mm}
  \item If
\begin{equation}\label{condS}
  \barDt_{f,j}(x_k,d_{k,j}) > \left(\frac{\varsigma\epsilon_j}{1+\nu}\right) \frac{\ord_k^j}{j!},
\end{equation}
  go to Step~2 with $j_k=j$.
  \end{enumerate}
  \vspace*{-2mm}
  Set $j_k=q$.
  \vskip 5pt
 \item[Step~2: Step computation.]
  If $r_k =\delta_k$, set $s_k = d_{k,j_k}$ and $\barDt_{f,j}(x_k,s_k)
  = \barDt_{f,j}(x_k,d_{k,j_k})$.
  Otherwise, compute a step $s_k$ such that $\|s_k\| \le r_k$ and
  \beqn{trqne-decrease2}
  \barDt_{f,j}(x_k,s_k) \ge \barDt_{f,j}(x_k,d_{k,j_k}).
\eeqn
 \item[Step~3: Function decrease estimation.] Compute the estimate $\barf(x_k)-\barf(x_k+s_k)$ of
   $f(x_k)-f(x_k+s_k)$.
    \vskip 5pt
\item[Step~4: Test of acceptance.] Compute
\begin{equation}\label{trqne-rhok-def}
  \rho_k = \frac{\barf(x_k) - \barf(x_k+s_k)}
                {\barDt_{f,j}(x_k,s_k)}.
\end{equation}
  If $\rho_k \geq \eta$ \textit{(successful iteration)}, then set
  $x_{k+1} = x_k + s_k$; otherwise \textit{(unsuccessful iteration)} set $x_{k+1} = x_k$.
  \vskip 5pt
\item[Step~5: Trust-region radius update.]
  Set
 \[
  r_{k+1} = \left\{ \begin{array}{ll}
  {}\frac{1}{\gamma}r_k, & \tim{if} \rho_k < \eta,\\
  {} \min[r_{\max},\gamma r_k], & \tim{if} \rho_k \ge  \eta,
  \end{array}\right.
  \]
  Increment $k$ by one and go to Step~1.
  \end{description}
}
  
\noindent
A feature of the \algn algorithm is that it uses an
adaptive strategy (in Step~1) to choose the model's degree in view of the
desired accuracy and optimality order. Indeed, the model of the objective
function used to compute the step is $\bart_{f,j_k}(x_k,s)$, whose degree $j_k$
can vary from an iteration to the other, depending on the ``order of (inexact)
optimality'' achieved at $x_k$ (as determined by Step~1). Also observe that,
if the trust-region radius is small (that is $r_k\le \theta$), the
optimality displacement $d_{k,j_k}$ is an approximate global minimizer of
the model within the trust region, which justifies the choice
$s_k=d_{k,j_k}$ in this case. If $r_k> \theta$, the step computation is
allowed to be fairly approximate as the only requirement for a step in the
trust region is \req{trqne-decrease2}.  This can be interpreted as a
generalization of the familiar notions of ``Cauchy'' and ``eigen'' points
(see \cite[Chapter~6]{ConnGoulToin00}).  In addition, note that, while nothing
guarantees that $f(x_k) \ge f(x_{k+1})$, the mechanism of the
algorithm ensures that $\barf(x_k) \ge \barf(x_{k+1})$.

The \al{TR$q$NE} algorithm generates a random process. Randomness occurs
because of the random noise present in the Taylor decreases and objective
function values, the former resulting itself from the randomly perturbed
derivatives values and, as the algorithm proceeds, from the random
realizations of the iterates $x_k$ and steps $s_k$. In the following analysis,
uppercase letters denote  random quantities, while lowercase ones denote realizations of these random quantities. Thus, 
given  $\omega\in \Omega$,   $x_k = X_k\rlz$, $g_k=G_k\rlz$, etc. 
In particular, we distinguish
\begin{itemize}
\item $\Dt_{f,j}(x,s)$, the value at a (deterministic) $x,s$ of the exact Taylor
decrement, that is of the Taylor decrement using the exact values of its derivatives at $x$;

\item $\barDt_{f,j}(x,s)=\barDt_{f,j}(x,s,\xi)$, the value at a (deterministic) $x,s$ of an inexact Taylor
decrement, that is of a Taylor decrement using the inexact values of 
its derivatives (at $x$) resulting from the realization of random
noise;

\item $\Dt_{f,j}(X,S)$, the random variable corresponding to the exact Taylor
  decrement taken at the random variables $X,S$;
  
\item$\barDT_{f,j}(X,S)$, the random variable giving the value of the Taylor
decrement using randomly perturbed derivatives, taken at
the random variables $X,S$.
\end{itemize}
Analogously, $F_k^0\eqdef F(X_k)$ and $F_k^s\eqdef F(X_k+S_k)$ denote the random variables associated with
the estimates of $f(X_k)$ and $f(X_k+S_k)$, with their realizations $f_k^0=\bar
f(x_k)=\bar f(x_k,\xi)$ and $f_k^s=\bar f(x_k+s_k)=\bar f(x_k+s_k,\xi)$.
Similarly, the iterates $X_k$, as well as the trust-region radiuses $R_k$, the
indeces $J_k$, the
optimality radiuses $\ORD_k$, displacements $D_{k,j}$
and the steps $S_k$  are random variables while
$x_k$, $r_k$, $j_k$, $\ord_k$,  $d_{k,j}$, 
and $s_k$ denote their realizations.
Hence, the \algn algorithm generates the random process
\begin{equation}
\{X_k, R_k, M_k, J_k, \ORD_k,
\{D_{k,j}\}_{j=1}^{J_k},
S_k, F_k \}\label{sprocess} 
\end{equation}
in which $X_0=x_0$ (the initial guess) and $R_0=r_0$ (the initial trust-region
radius) are deterministic quantities, and where
  \[
M_k = \{\overline{\nabla_x^1 f}(X_k), \ldots, \overline{\nabla_x^{j_k} f}(X_k)\}
\mand
F_k = \{F(X_k), F(X_k+S_k)\}.
\]

\subsection{The probabilistic setting}

We now state our probabilistic assumptions on the \algn algorithm. For $k\ge
0$, our assumption on the past is formalized by considering $\sigax$ the
$\sigma$-algebra induced by the random variables $M_0$, $M_1$,..., $M_{k-1}$
and $F_0^0$, $F_0^s$, $F_1^0$, $F_1^s$, ..., $F_{k-1}^0$, $F_{k-1}^s$ and let
$\calA_{k-1/2}$ be that induced by $M_0$, $M_1$,..., $M_k$ and $F_0^0$,
$F_0^s$, ..., $F_{k-1}^0$, $F_{k-1}^s$, with $\calA_{-1}=\sigma(x_0)$.

We first define an event ensuring that the model is accurate enough at
iteration $k$. At the end of Step 2 of this iteration and given
$J_k\in\ii{q}$, we now define,
\begin{align}
  \calM_{k,j}^{(1)}&=\left\{
  \phi_{f,j}^{\ORD_k}(X_k)
  \le \left(\frac{1+\nu}{\varsigma}\right) \barDT_{f,j}(X_k,D_{k,j})
                     \right\} \ms (j\in\ii{J_k}),\nonumber\neol
  \calM_k^{(2)}&=\left\{
  (1-\nu) \barDT_{f,J_k}(X_k,S_k) \le
  \Dt_{f,J_k}(X_k,S_k)
  \le (1+\nu) \barDT_{f,J_k}(X_k,S_k)
                     \right\},\nonumber\neol
\calM_k\;&=\left(\bigcap_{j\in\ii{J_k}}\calM_{k,j}^{(1)} \right)\cap \calM_k^{(2)}.\label{Mk}
\end{align}
The event $\calM_{k,j}^{(1)}$ occurs when the $j$-th order optimality measure
($j\le j_k$) at iteration $k$ is meaningful, while $\calM_k^{(2)}$ occurs when
this is the case for the model decrease.  At first sight, these events may
seem only vaguely related to the accuracy of the function's derivatives but a
closer examination gives the following sufficient condition for $\calM_k$ to
happen.

\llem{suff-acc-conds}{At iteration $k$ of any realization, the inequalities
  defining the event $\calM_k$ are satisfied if,
  for $j \in \ii{j_k}$ and $\ell \in \ii{j_k}$
  \beqn{the-cond-0}
  \hspace*{-21mm}
  \|\left(\overline{\nabla_x^\ell f}(x_k) -\nabla_x^\ell f(x_k)\right)[s_k]^\ell\|
  \leq \frac{\nu}{2} \barDt_{f,j_k}(x_k,s_k)
  \eeqn
  and
  \beqn{the-cond-2}
  \|\left(\overline{\nabla_x^\ell f}(x_k) -\nabla_x^\ell f(x_k)\right)[\hat{d}_{k,j}]^\ell\|
  \leq 
	{\frac{\nu}{2}  \barDt_{f,j}(x_k,\hat{d}_{k,j}),}
  \eeqn
  where  
   \beqn{hatd}
  \hat{d}_{k,j} = \argmax_{\|d\| \le \ord_k} \Dt_{f,j}(x_k,d).
  \eeqn
}

\proof{
If \req{the-cond-2} holds, we have that, for every $j\in \ii{j_k}$ and
$v_k\in\{\hat{d}_{k,1},\ldots,\hat{d}_{k,j}\}$,  
with  $\hat{d}_{k,\ell}$, $\ell=1,\ldots,j$, given in \eqref{hatd},
\beqn{rel-err}
\begin{array}{lcl}
|\barDt_{f,j}(x_k,v_k)-\Dt_{f,j}(x_k,v_k)|
& \le & \bigsum_{\ell=1}^j \bigfrac{1}{\ell!}\,
         \left\|\left(\overline{\nabla_x^\ell f}(x_k) -\nabla_x^\ell f(x_k)\right)[v_k]^\ell\right\|\\*[2ex]
& \le &\bigfrac{1}{2}\nu  \barDt_{f,j}(x_k,v_k) \bigsum_{\ell=1}^{j} \frac{1}{\ell!}\,\\*[2ex]
& <  & \nu \,  \barDt_{f,j}(x_k,v_k) 
\end{array}
\eeqn
where we have used the bound $\bigsum_{\ell=1}^{j} \frac{1}{\ell!} < e-1 <2$.
Now note that the definition of $\barphi_{f,j}^{\ord_k}(x_k)$
in \req{phi-def-bar}, \req{rel-err} for $v_k = \hat{d}_{k,j}$
and \eqref{globmax_inexact} imply that, for any $j\in\ii{j_k}$,
\[
\begin{array}{lcl}
\phi_{f,j}^{\ord_k}(x_k)=\Dt_{f,j}(x_k,\hat{d}_{k,j})
& \leq & \barDt_{f,j}(x_k,\hat{d}_{k,j})+|\Dt_{f,j}(x_k,\hat{d}_{k,j})-\barDt_{f,j}(x_k,\hat{d}_{k,j})|\\*[2ex]
& \leq & \big(1+\nu\big)\barDt_{f,j}(x_k,\hat{d}_{k,j})\\*[2ex]
& \leq & \big(1+\nu\big) \max_{\|d\|\le\ord_k}\barDt_{f,j}(x_k,d)\\*[2ex]
& = & \big(1+\nu\big)\barphi_{f,j}^{\ord_k}(x_k)\\*[2ex]
& \le &  \left(\bigfrac{1+\nu}{\varsigma}\right) \barDt_{f,j}(x_k,d_{k,j}).
\end{array}
\]
Hence the 
inequality in the definition of $\calM_{k,j}^{(1)}$ holds for
$j\in \ii{j_k}$. The proof of the inequalities defining $\calM_{k}^{(2)}$ is
analog to that of \req{rel-err}. We have from \req{the-cond-0} that
\beqn{rel-err-2}
\begin{array}{lcl}
|\barDt_{f,j_k}(x_k,s_k)-\Dt_{f,j_k}(x_k,s_k)|
& \le & \bigsum_{\ell=1}^{j_k} \bigfrac{1}{\ell!}\,
         \left\|\left(\overline{\nabla_x^\ell f}(x_k) -\nabla_x^\ell f(x_k)\right)[s_k]^\ell\right\|\\*[2ex]
& \le &\bigfrac{1}{2}\nu  \barDt_{f,j_k}(x_k,s_k) \bigsum_{\ell=1}^{j_k} \frac{1}{\ell!}\,\\*[2ex]
& <  & \nu \,  \barDt_{f,j_k}(x_k,s_k) 
\end{array}
\eeqn
where we have again used the bound $\bigsum_{\ell=1}^{j_k} \frac{1}{\ell!} < 2$.
} 

\noindent
This result immediately suggests a few comments.
\begin{itemize}
\item The conditions \req{the-cond-0}-\req{the-cond-2} are merely sufficient,
  not necessary.  In particular, they ignore any possible cancellation of errors
  between terms of the Taylor expansion of different degree.
\item We note that \req{the-cond-0}-\req{the-cond-2} require the $\ell$-th
  derivative to be relatively accurate along a finite and limited set of
  directions, independent of problem dimension.
\item Since $\|d_{k,j}\|$ and $\|\hat{d}_{k,j}\|$ are bounded by $\ord_k \le
  \theta \le 1$, we also note that the accuracy required by these conditions
  decreases when the degree $\ell$ increases.  Moreover, for a fixed degree,
  the request is weaker for small displacements (a typical situation when a
  solution is approached) than for large ones.
\item Requiring
  \beqn{the-cond-2-n}
  \|\overline{\nabla_x^\ell f}(x_k) -\nabla_x^\ell f(x_k)\|
  \leq \frac{\nu}{2\|\hat{d}_{k,j}\|^\ell} {\barDt_{f,j}(x_k,\hat{d}_{k,j}),}
  \eeqn
  {instead of \req{the-cond-2}} is of course again \emph{sufficient} to ensure the desired conclusions.
  These conditions are reminiscent of the conditions required in
  \cite{BlanCartMeniSche19} for the STORM algorithm with $p=2$, namely that,
  for some constant $\kappa_\ell$ and all $y$ in the trust-region
  $\{ y \in \Re^n \mid \|y-x_k\| \le r_k\}$, 
$$
  \|\overline{\nabla_x^\ell f}(y) -\nabla_x^\ell f(y)\|\le \kappa_\ell \,\, r_k^{3-\ell}
  \ms (\ell\in\{0,1,2\}).
 $$
  This latter  condition  is however stronger 
	than \req{the-cond-0}--\req{the-cond-2}
  because it insists on a uniform accuracy guarantee in the full-dimensional
  trust region.
\end{itemize}

\noindent
Having considered the accuracy of the model, we now turn to that on the
objective function's values. At the end of Step $3$ of the $k$-th iteration,
we define the event
\beqn{Fk}
\calF_k=\left\{|\Df(X_k,S_k)-\DF(X_k,S_k)|\le 2\nu \barDT_{f,j_k}(X_k,S_k)\right\}
\eeqn
where $\Df(X_k,S_k) \eqdef f(X_k)-f(X_k+S_k)$ and $\DF(X_k,S_k) \eqdef
F(X_k)-F(X_k+S_k)$. This occurs when the difference in function values used in
the course of iteration $k$ are reasonably accurate relative to the model decrease
obtained in that iteration. Note that, because of the triangular inequality,
\begin{align*}
|\Df(X_k,S_k)-\DF(X_k,S_k)|
&  =  |(f(X_k)-f(X_k+S_k)) - (F(X_k)-F(X_k+S_k))|\\
& \le |f(X_k)-F(X_k)| + |f(X_k+S_k)-F(X_k+S_k)|
\end{align*}
so that $\calF_k$ must occur if both terms on the right-hand side are bounded
above by $\nu\barDT_{f,j_k}(X_k,S_k)$.
Combining accuracy requests on model and function values, we define
\beqn{E-def}
\calE_k\eqdef \calF_k \cap \calM_k
\eeqn
and say that iteration $k$ is \textit{accurate} if
$\indicE =\indicF\indicM = 1$ and the iteration $k$ is \textit{inaccurate} if
$\indicE = 0$.
Moreover, we say that the iteration $k$ has accurate model if $\indicM=1$
and that iteration $k$ has accurate function estimates if $\indicF=1$.  
Finally we let
\[
p_{\calM_k} \eqdef \prob\Big[\calM_k \mid \sigax\Big], \ \  p_{\calF_k} \eqdef \prob\Big[\calF_k \mid \calA_{k-1}\Big].
\]

We will verify in what follows that the \algn algorithm does progress towards
an approximate minimizer satisfying \req{approx-min} as long as the following holds.

\begin{description}
\item[AS.3] There exists $\alpha_*, \gamma_*\in (\half,1]$ such that $p_*=\alpha_*\gamma_*>\half$,
\beqn{prob-ass}
p_{\calM_k}\ge \alpha_*,
\ms
p_{\calF_k} \ge \gamma_*
\mand
\expect\Big[\indicS(1-\indicF)\Df(X_k,S_k) \mid \calA_{k-1} \Big] \ge 0,
\eeqn
where $\calS_k$ is the event $\calS_k\eqdef \{  \mbox{iteration $k$ is successful}\}$.
\end{description}

\noindent
We notice that due to the tower property for conditional expectations
$$
\prob\Big[\calF_k \mid \sigax\Big]=\expect \Big[\indicF\mid \sigax\Big]=\expect \Big[\expect \Big[\indicF\mid  \calA_{k-\half}\Big]\mid \sigax \Big].
 $$
and hence that assuming, as in \cite{PaquSche20} and \cite{BlanCartMeniSche19},
\[
\prob\Big[\calF_k \mid \calA_{k-\frac{1}{2}}\Big] > \gamma_*
\]
is stronger than assuming  $p_{\calF_k} \ge \gamma_*$.
Similarly, 
\beqn{expectstronger}
\expect\Big[\indicS(1-\indicF)\Df(X_k,S_k) \mid  \calA_{k-\frac{1}{2}}\Big] \ge 0
\tim{implies}
\expect\Big[\indicS(1-\indicF)\Df(X_k,S_k) \mid  \calA_{k-1}\Big] \ge 0.
\eeqn

Assuming AS.3 is not unreasonable as it merely requires that an accurate model  {and accurate functions} 
``happen more often than not'', and that the discrepancy between true and
inexact function values at successful iterations does not, on average, prevent
decrease of the objective function.  If either of these
condition fails, it is easy to imagine that the \algn algorithm could be completely
hampered by noise and/or diverge completely.  Because the last condition in
\req{prob-ass} is less intuitive, we now show that it can be realistic in the
specific context where reasonable assumptions are made on the (possibly
extended) cumulative distribution of the error on the function decreases
(conditioned to $\calA_{k-\half}$).

\lthm{expect-ass-th}{
Let $G_k: \Re_+\rightarrow [0,1]$ be a differentiable monotone
increasing random function which is measurable for $\calA_{k-1}$ and such that   
\begin{eqnarray}
& & G_k(0)=0 \tim {and} \lim_{\tau\rightarrow \infty} G_k(\tau)=1,\label{g0}\\
& & \lim_{\tau\rightarrow \infty} \tau\,(1-G_k(\tau) )=0\label{g1},\\
& & \int_{0}^\infty (1-G_k(\tau))\,d\tau<\infty, \label{g2}
\end{eqnarray}
and
\beqn{pb3}
\prob\Big[ \DF(X_k,S_k)-\Df(X_k,S_k) >  \tau  \mid \calA_{k-\half}\Big ] \le 1-G_k(\tau)
\eeqn
for $\tau > 0$. Then, 
\begin{equation}\label{expectation}
\expect \Big[ \indicS(1-\indicF)(f(X_k)-f(X_{k+1})) \mid \calA_{k-1} \Big] \ge 0
\end{equation}
for each $k$ such that
\[
\barDT_{f,j_k}(X_k,S_k) \ge \frac{1}{\eta}\int_{0}^\infty (1-G_k(\tau))\,d\tau.
\]
}

\proof{ 
Consider $\omega\in \Omega$, an arbitrary realization of the stochastic process defined
by the \al{TR$q$NE} algorithm.
Suppose first that  $\expect\Big[\indicS(1-\indicF) \mid \calA_{k-\half}\Big]\rlz = 0$.
We then deduce that
\beqn{case0}
\expect\Big[\indicS(1-\indicF)(f(X_k)-f(X_{k+1})) \mid \calA_{k-\half} \Big]\rlz =0.
\eeqn
Assume therefore that
\begin{equation} \label{bound-indicF-ass}
\expect \Big[\indicS(1-\indicF) \mid \calA_{k-\half}\Big]\rlz= \bar p_k
\end{equation}
for some $\bar p_k>0$. To further simplify notations, set
\beqn{Ek}
\barDT_k\eqdef\barDT_{f,j_k}(X_k,S_k)\eta
\tim{and}
E_k =\DF(X_k,S_k)-\Df(X_k,S_k).
\eeqn
If we define $\calI \eqdef \{ E_k\rlz> 0\}$,
the definition of successful iterations, \req{trqne-rhok-def} and the
triangular inequality then imply that,  if
$\indic{\calS_k\rlz}=1$ then
\beqn{Df-vs-DT}
\Df(X_k,S_k)\rlz = \DF(X_k,S_k)\rlz - E_k\rlz
\ge \eta \barDT_k\rlz -  \indicI E_k\rlz.
\eeqn
This in turn ensures that
\begin{align}\label{expect-bound}
  \expect\Big[\indicS(1-\indicF)\Df(X_k,S_k) \mid \calA_{k-\half} \Big]\rlz
&\ge \eta\barDT_k\rlz \, \expect \Big[\indicS(1-\indicF) \mid \calA_{k-\half}\Big]\rlz \neol
& - \expect \Big[\indicS(1-\indicF) \indic{{\cal I}}E_k\,\mid \calA_{k-\half}\Big]\rlz.
\end{align}
Moreover, we have that
\begin{align*}
\expect \Big[\indicS&(1-\indicF) \indicI E_k\,\mid \calA_{k-\half}\Big]\rlz\\
& = \expect \Big[\indicS(1-\indicF) \indicI E_k\,
  \mid \calA_{k-\half}, \calS_k\cap\no{\calF_k}\Big]\rlz\,
   \cdot \prob[\calS_k\cap\no{\calF_k}\,\mid \calA_{k-\half}]\rlz\\
& \hspace*{1cm}+\expect \Big[\indicS(1-\indicF) \indicI E_k\,
  \mid \calA_{k-\half}, \no{(\calS_k\cap\no{\calF_k})}\Big]\rlz\,
   \cdot \prob[\no{(\calS_k\cap\no{\calF_k})}\,\mid \calA_{k-\half}]\rlz\\
& = \bar p_k \, \expect \Big[\indicI E_k\,\mid \calA_{k-\half}, \calS_k\cap\no{\calF_k}\Big]\rlz\\
& \le \bar p_k\, \expect \Big[\indicI E_k\,\mid \calA_{k-\half}\Big]\rlz,
\end{align*}
where we used the fact that $[\indicS(1-\indicF)]\rlz=0$ whenever
${(\no{\calS_k}\cup{\calF_k})}\rlz$ happens, \req{bound-indicF-ass} to derive the second equality, and the bound
$\indicI E_k\rlz \ge 0$ to obtain the final inequality. 
Now, \eqref{pb3} implies that, for $\tau>0$
$$
\prob\Big[ \indicI E_k >  \tau \,\mid \calA_{k-\half}\Big]\rlz
=\prob\Big[ E_k >\tau\,\mid \calA_{k-\half} \Big]\rlz
\le 1-g_k(\tau)
$$
where $g_k(\tau) \eqdef G_k\rlz(\tau)$,
and thus
\[
\prob\Big[\indicS(1-\indicF)\indicI E_k>  \tau \,\mid \calA_{k-\half}\Big ]\rlz
\le (1-g_k(\tau))\,\bar p_k
= \bar p_k \int_{\tau}^\infty g_k'(t)dt.
\]
Then, employing \req{g0}--\req{pb3},  and integrating by parts
\[
\expect \Big[\indicS(1-\indicF)\indicI E_k \,\mid \calA_{k-\half}\Big]\rlz
\le \bar p_k \int_{0}^\infty t\, g_k'(t)dt
=\bar p_k \int_{0}^\infty (1-g_k(t))dt
< \infty.
\]
Finally, using \req{expect-bound},
\begin{eqnarray*}
  \expect\Big[\indicS(1-\indicF)\Df(X_k,S_k) \mid \calA_{k-\half} \Big]\rlz
&\ge \bar p_k \left[ \eta\barDT_k\rlz -  \bigint_{0}^\infty (1-g_k(t))\,dt \right]
\end{eqnarray*}
and thus
\begin{equation}\label{expectation-rlz}
\expect \Big[ \indicS(1-\indicF)(f(X_k)-f(X_{k+1})) \mid \calA_{k-\half} \Big]\rlz \ge 0
\end{equation}
holds when
\[
\barDT_k\rlz\ge \frac{1}{\eta}\int_{0}^\infty (1-g_k(t))\,dt.
\]
Combining \req{case0} and \req{expectation-rlz} and taking into account that
$\omega$ is arbitrary give that
\[
\expect \Big[ \indicS(1-\indicF)(f(X_k)-f(X_{k+1})) \mid \calA_{k-\half} \Big]
\ge 0,
\]
which, in view of \req{expectstronger}, yields \req{expectation}.
}

\noindent
Note that the assumptions of the theorem are for instance satisfied for
the exponential case where $G_k(\tau)= e^{-T \tau}$ for $T >0$ and measurable
for $\calA_{k-1}$. We will return to this result in Section~\ref{noise-s} and
discuss there the condition that $\barDT_{f,j_k}(X_k,S_k)$ should be
sufficiently large.

\numsection{Worst-case evaluation complexity}\label{section:analysis}  

We now turn to the evaluation complexity analysis for the \algn algorithm, whose aim is to derive
a bound on the expected number of iterations for which optimality fails. This
number is given by
\beqn{Neps-def}
N_\epsilon
\eqdef \inf \left\{
k\geq 0 ~\mid~  \phi_{f,j}^{\ORD_{k,j}}(X_k) \leq \epsilon_j
\frac{\ORD_{k,j}^j}{j!} \tim{for} j \in \ii{q} 
\right\}.
\eeqn

\noindent
We first state a crucial lower bound on the model decrease, in the spirit of 
\cite[Lemma~3.4]{CartGoulToin20c}.

\llem{trqne-phihat-s}{Consider any realization of the algorithm and assume that $\calM_k$ occurs. 
Assume that \req{approx-min} fails at iteration $k$.
Then, {there exists a $j_k\in \{1, \ldots, q\}$ such that
$ \barDt_{f,j_k}(x_k,d_{k,j_k}) > \varsigma \epsilon_{j_k} \ord_k^{j_k} / (j_k!(1+\nu))$}
at Step~1 of the iteration. Moreover,
{
\beqn{unco-trqne-DTjk}
\barDt_{f,j_k}(x_k,s_k) \geq \sphi_{f,k}\frac{\ord_k^{j_k}}{j_k!}
\eeqn
where
\beqn{trqne-phihat-def}
\widehat{\phi}_{f,k}
\eqdef \bigfrac{j_k!\,\,\barDt_{f,j_k}(x_k,d_{k,j_k})}{\ord_k^{j_k}}
> \frac{\varsigma\epsmin}{1+\nu}.
\eeqn
}
}

\proof{ 
We proceed by contradiction and assume that 
\begin{equation}\label{assurdo}
\barDt_{f,j}(x_k,d_{k,j})\le \frac{\varsigma\epsilon_j}{1+\nu} \,\frac{\ord_k^j}{j!},
\end{equation}
for all $j\in\ii{q}$. Since $\calM_k$ occurs, we have that, for all $j\in\ii{q}$,
\[
 \phi_{f,j}^{\ord_k}(x_k)
 \le \left(\bigfrac{1+\nu}{\varsigma}\right) \barDt_{f,j}(x_k,d_{k,j})
 \le \epsilon_j \,\frac{\ord_k^j}{j!},\quad j\in\{1,...,q\},
\]
which contradicts the assumption that \eqref{approx-min} does not hold for
$x_k$ and $\ord_k$. The bound \eqref{unco-trqne-DTjk}
directly results from
\[
\barDt_{f,j_k}(x_k,s_k)\ge \barDt_{f,j_k}(x_k,d_{k,j_k})
= \sphi_{f,k}\frac{\ord_k^{j_k}}{j_k!},
\]
where we have used \req{trqne-decrease2} to derive the first inequality
and the definition \eqref{trqne-phihat-def} to obtain the equality. 
The
rightmost inequality in \req{trqne-phihat-def} trivially follows from the negation of \req{assurdo} and
\req{epsmin-omega-def}.
} 

\noindent
We now search for conditions ensuring that the iteration is successful.
For simplicity of notation, given $\vartheta_{f,j}$, $j\in \{1,\ldots, q\}$, as in \req{f-holder},  
we define
\beqn{unco-Lf-def}
\vartheta_f \eqdef \max[1, \max_{j\in\ii{q}} \vartheta_{f,j}]. 
\eeqn

\llem{trqne-condsucc-l}{
Suppose that AS.1 holds. Consider any realization of the algorithm and
suppose that \eqref{approx-min} does not hold for $x_k$ and $\ord_k$ and that
$\calE_k$ occurs. If 
\beqn{trqne-condsucc}
r_k
\leq \rr
\eqdef\min\left\{\theta,\frac{\varsigma(1-\eta)}
                  {4(1+\nu)\vartheta_f}\,\epsilon_{\min} \right\}
= \frac{\varsigma(1-\eta)}{4(1+\nu)\vartheta_f}\,\epsilon_{\min}
\eqdef \kappa_r \epsmin,
\eeqn
$\kappa_r \in (0,1),$
holds,
then iteration $k$ is successful. 
}

\proof{
First, note that the minimum in \eqref{trqne-condsucc} is attained at $\kappa_r
\epsmin$ since $\theta \ge \epsmin$ 
and $\kappa_r\in
(0,1)$. Suppose now that \eqref{trqne-condsucc} holds,  which  implies that $\ord_k =
\min[\theta,r_k] = r_k$. Let $j_k$ be as in Lemma \ref{trqne-phihat-s}, and
denote  $\Df(x_k,s_k)= f(x_k)-f(x_k+s_k)$, $\barDf(x_k,s_k)=\barf(x_k) - \barf(x_k+s_k)$.

Using \req{trqne-rhok-def}, the triangle inequality  and $\indicE=\indicM\indicF=1$,
we obtain
\begin{align*}
|\rho_k - 1|
& = \left|\frac{\barDf(x_k,s_k)-\barDt_{f,j_k}(x_k,s_k)}{\barDt_{f,j_k}(x_k,s_k)}\right|\\
& \le \frac{\left|\Df(x_k,s_k)-\Dt_{f,j_k}(x_k,s_k)\right|}{\barDt_{f,j_k}(x_k,s_k)}
                + \frac{\left|\Dt_{f,j_k}(x_k,s_k)-\barDt_{f,j_k}(x_k,s_k)\right|}{\barDt_{f,j_k}(x_k,s_k)}\\
& \hspace*{7mm} +\frac{\left|\barDf(x_k,s_k)-\Df(x_k,s_k)\right|}{\barDt_{f,j_k}(x_k,s_k)}\\
& {\le} \frac{|f(x_k+s_k)-t_{f,j_k}(x_k,s_k)|}{\barDt_{f,j_k}(x_k,s_k)}+ \frac{3\nu\barDt_{f,j_k}(x_k,s_k)}{\barDt_{f,j_k}(x_k,s_k)}.
\end{align*}

\noindent
Invoking \req{tech-resf}, the bound $\|s_k\|\leq r_k=\ord_k$,
\req{unco-Lf-def}, \req{unco-trqne-DTjk} and $\nu \leq \frac 1 4 (1-\eta)$ we get
$$
|\rho_k - 1|
< \bigfrac{\vartheta_fr_k}{\sphi_{f,k}} +\frac{3}{4}(1-\eta).
$$
Using  \eqref{trqne-phihat-def} and \req{trqne-condsucc}  we deduce that
{
\begin{equation}
\label{chain}
|\rho_k - 1|
 \le  \bigfrac{(1+\nu)\vartheta_fr_k}{\varsigma\epsilon_{\min}} +\frac{3}{4}(1-\eta)
 \le 1-\eta.
 \end{equation}
}
Thus, $\rho_k \geq \eta$ and the iteration $k$ is successful.
} 

\noindent
The following crucial lower bound on $\barDT_{f,j_k}(x_k,s_k)$ for
accurate iterations $k$ can now be proved.

\llem{trqne-model-decrease-2}{
{Suppose that AS.1 holds. Consider any realization of the algorithm and
  suppose that \eqref{approx-min} does not hold for $x_k$ and $\ord_k$,  
  $\calE_k$ occurs, and  that $r_k\ge \bar r$ with $\bar r_k$ defined in \req{trqne-condsucc}. Then
  \beqn{trqne-model-decrease2}
\barDt_{f,j_k}(x_k,s_k)=\bart_{f,j_k}(x_k,0)-\bart_{f,j_k}(x_k,s_k)
{> \frac{\varsigma}{q!} \left(\kappa_\ord\epsmin\right)^{q+1}},
\eeqn
where $\kappa_\ord \in (0,1)$ is defined by
\beqn{trqne-kappddef2}
\kappa_\ord \eqdef \frac{\kappa_r}{1+\nu}
\eeqn
with $\kappa_r$ defined in (\ref{trqne-condsucc}).
} }

\proof{
Let $j_k$ be as in Lemma \ref{trqne-phihat-s}. By 
\eqref{unco-trqne-DTjk}, \eqref{trqne-phihat-def} we obtain
\[
\barDt_{f,j_k}(x_k,s_k) {> \frac{\varsigma\epsmin}{1+\nu} \,\frac{\ord_k^q}{q!}.}
\]
If $r_k > \theta$  then $\ord_k=\theta$ and the bound $\theta \geq \epsmin$ implies
\[
\barDt_{f,j}(x_k,s_k) { > \frac{\varsigma\epsmin^{q+1}}{q!(1+\nu)}.}
\]
Thus \req{trqne-model-decrease2} holds by definition
of $\kappa_\ord$ and the fact that $\kappa_r\in (0,1)$.
If $\bar r <r_k \leq \theta$, then $\ord_k= r_k$.
The proof is completed by noting that the form of $\bar r$ in (\ref{trqne-condsucc})
gives that $r_k>  \kappa_r \epsmin$.\\
} 

\subsection{Bounding the expected number of steps with
  \boldmath$R_k \le \rr$}

We now return to the general stochastic process generated by the \algn algorithm
and  bound  the expected number of steps in $N_{\epsilon}$ from above. For this
purpose, let us define, for all $0 \le k 
\le N_{\epsilon}-1$,
the events 
\[
\Lambda_k \eqdef \{ R_k > \rr\},
\qquad
\no{\Lambda}_k \eqdef \ \{R_k \le \rr\},
\]
where $\rr$ is given by \req{trqne-condsucc}, and let 
\beqn{defnsigma}
N_\Lambda \eqdef  \sum_{k=0}^{N_\epsilon-1}\indicL,
\qquad
N_{\no{\Lambda}}
\eqdef \sum_{k=0}^{N_\epsilon-1}\indicLc,
\eeqn
be the number of steps, in the stochastic process induced by the \algn
algorithm and before
$N_{\epsilon}$, such that $R_k> \rr$ or $R_k\le\rr$,
respectively. In what follows we suppose that AS.1--AS.2 hold.

\noindent
An upper bound on  $\expect\big[N_{\no{\Lambda}}\big]$ can be derived as follows.
\begin{itemize}
\item[(i)] We apply \cite[Lemma~2.2]{CartSche17} to deduce that, for any
  $\ell\in\iibe{0}{N_\epsilon-1}$ and for all realizations of Algorithm \ref{algoS}, one has that
\beqn{CSlemma22}
\sum_{k=0}^\ell \indicLc\indicS \le \frac{\ell+1}{2}.
\eeqn
\item[(ii)] Both
  $\sigma(\indicL)$ and $\sigma(\indicLc)$
  belong to $\sigax$, because the random variable $\Lambda_k$ is fully determined
  by the first $k-1$ iterations of the \algn algorithm. Setting $\ell=N_{\epsilon}-1$  we can rely on
  \cite[Lemma $2.1$]{CartSche17} (with $W_k=\indicLc$), whose proof is detailed in the appendix, to deduce that  
\beqn{CSlemma21}
\expect\left[ \sum_{k=0}^{N_\epsilon-1}\indicLc\indicE \right]
\geq \,\pM\,\expect\left[ \sum_{k=0}^{N_\epsilon-1}  \indicLc \right].
\end{equation}
\item[(iii)] As a consequence, given that Lemma \ref{trqne-condsucc-l}
ensures that each iteration $k$ where $\calE_k $ occurs and
$r_k\le \rr$ is successful, we have that 
$$
\sum_{k=0}^{N_\epsilon-1} \indicLc \indicE 
\leq  \sum_{k=0}^{N_\epsilon-1} \indicLc \indicS
\leq \frac{N_\epsilon}{2},
$$
in which the last inequality follows from \eqref{CSlemma22}, with
$\ell=N_\epsilon-1$. Taking expectation in the above inequality, using
\eqref{CSlemma21} and recalling the rightmost definition in \eqref{defnsigma},
we obtain, as in \cite[Lemma~2.3]{CartSche17}, that
\beqn{bound_Ns}
\expect[N_{\no{\Lambda}}] \leq \frac{1}{2\pM}\expect[N_\epsilon].
\eeqn
\end{itemize} 
\noindent

\subsection{Bounding the expected number of steps with
   \boldmath$R_k > \rr$}
For analyzing $\expect[N_\Lambda]$, where $N_\Lambda$ is defined in
\req{defnsigma}, we now introduce the following variables. 
\vspace{0.1cm}
\begin{definition}
Consider the random process \eqref{sprocess} generated by
the \al{TR$q$NE} algorithm and define: 
\vspace{0.1cm}
\begin{equation}
\begin{array}{lcl}
  \bullet\ \overline{\Lambda}_k &\hspace*{-3mm}=&
     \{\tim{iteration $k$ is such that $R_k\ge \rr$}\};\\*[1.5ex]
\bullet\ N_I &\hspace*{-3mm}=& \bigsum_{k=0}^{N_{\epsilon}-1}\indicbL \indicEc:
       \mbox{the number of inaccurate iterations with $R_k\ge \rr$};\\*[1.5ex]
\bullet\ N_A &\hspace*{-3mm}=& \bigsum_{k=0}^{N_{\epsilon}-1}\indicbL \indicE:
       \mbox{the number of accurate iterations with $R_k\ge \rr$};\\*[1.5ex]
\bullet\ N_{AS} &\hspace*{-3mm}=& \bigsum_{k=0}^{N_{\epsilon}-1}\indicbL\indicE\indicS:
        \mbox{the number of accurate successful iterations with $R_k\ge \rr$};\\*[1.5ex]
\bullet\ N_{AU} &\hspace*{-3mm}=& \bigsum_{k=0}^{N_{\epsilon}-1}\indicL\indicE\indicSc:
        \mbox{the number of accurate unsuccessful iterations with $R_k> \rr$};\\*[1.5ex]
\bullet\ N_{IS} &\hspace*{-3mm}=& \bigsum_{k=0}^{N_{\epsilon}-1}\indicbL\indicEc\indicS:
        \mbox{the number of inaccurate successful iterations with $R_k\ge \rr$};\\*[1.5ex]
\bullet\ N_S &\hspace*{-3mm}=& \bigsum_{k=0}^{N_{\epsilon}-1}\indicbL\indicS:
        \mbox{the number of successful iterations with $R_k\ge \rr$};\\*[1.5ex]
\bullet\ N_U &\hspace*{-3mm}=& \bigsum_{k=0}^{N_{\epsilon}-1}\indicL\indicSc:
        \mbox{the number of unsuccessful iterations with $R_k > \rr$}.
\end{array}
\end{equation}
\label{def2}
\end{definition}

\noindent
Observe that $\overline{\Lambda}_k$ is the ``closure'' of $\Lambda_k$ in
that the inequality in its definition is no longer strict.
\noindent
We immediately notice that an upper bound on $\expect[N_\Lambda]$
is available, once an upper bound on $\expect[N_I]+\expect[N_A]$ is
known, since  
\begin{equation}
\label{PlanEN}
\expect[N_\Lambda]
\leq \expect\left[\sum_{k=0}^{N_{\epsilon}-1}\indicbL \right]
= \expect\left[
   \sum_{k=0}^{N_{\epsilon}-1}\indicbL\indicEc
   +\sum_{k=0}^{N_{\epsilon}-1}\indicbL\indicE\right]
= \expect[N_I]+\expect[N_A].
\end{equation}
Using again \cite[Lemma $2.1$]{CartSche17} (with
$W_k=\indicbL$) to give an upper bound on
$\expect[N_I]$, we obtain the following result, whose proof is detailed in the appendix.

\llem{a1}{\cite[Lemma~2.6]{CartSche17}
Suppose that AS.1-AS.3 hold and let $N_I$, $N_A$ be defined as in Definition
\ref{def2} in the context of the stochastic process \eqref{sprocess}
generated by the \algn algorithm. Then,
\begin{equation}
\label{EM1}
\expect[N_I]\le \frac{1-\pM}{\pM} \,\expect[N_A].
\end{equation}
}

\noindent
Turning to  the upper bound for $\expect[N_A]$, we observe that 
\begin{equation}
\label{defEM2}
\expect[N_A]
  =  \expect[N_{AS}]+\expect[N_{AU}]
\leq \expect[N_{AS}]+ \expect[N_U].
\end{equation}
\noindent
Hence, bounding $\expect[N_A]$ can be achieved by providing upper bounds on
$\expect[N_{AS}]$ and $\expect[N_U]$. Regarding the latter, we first
note that the process induced by the \algn algorithm ensures that $R_k$
is increased by a factor $\gamma$ on successful steps and decreased by the
same factor on unsuccessful ones.
Consequently, based on \cite[Lemma $2.5$]{CartSche17}, we obtain the
following bound. 

\llem{a2}{
  For any $\ell\in\{0,...,N_{\epsilon}-1\}$ and for
  all realisations of Algorithm \ref{algoS}, we have that
\[
\sum_{k=0}^\ell \indicL\indicSc
\leq \sum_{k=0}^{\ell}\indicbL\indicS
  +\left\lceil\Big|\log_{\gamma^{-1}}\left( \frac{\rr}{r_0} \right)\Big|\right\rceil=\sum_{k=0}^{\ell}\indicbL\indicS
  +\left\lceil\log_{\gamma}\left( \frac{r_0}{\rr} \right)\right\rceil.
\]
}

\noindent
From the inequality stated in the previous lemma with $\ell=N_{\epsilon}-1$, recalling Definition
\ref{def2} and taking expectations, we therefore obtain that
\begin{equation}
\label{boundEU}
\expect[N_U]
\le \expect[N_S]+ \left\lceil\log_{\gamma}\left( \frac{r_0}{\rr} \right)\right\rceil
=\expect[N_{AS}]+\expect[N_{IS}]+\left\lceil\log_{\gamma}\left( \frac{r_0}{\rr} \right)\right\rceil.
\end{equation}

\noindent
An upper bound on $\expect[N_{AS}]$ is given by the following lemma.\\

\llem{LemmaNAS}{ Suppose that AS.1, AS.2  and AS.3 hold.
  Then we have that 
\begin{equation}
\label{boundEN1}
\expect[N_{AS}]
\le \frac{2q!(f_0-\flow)}{\varsigma\eta
  \left(\kappa_\delta\epsmin\right)^{q+1}} +1,
\end{equation}
where $\kappa_{\ord}$ is defined in \req{trqne-kappddef2}.
}

\proof{
Consider any realization of the \algn algorithm.
\begin{itemize}
\vskip 5pt
\item[\textit{i)}] If iteration $k$ is successful and the
  functions are accurate (i.e., $\indicS\indicF=1$)
  then {\eqref{trqne-rhok-def},   \req{epsmin-omega-def} and \req{Fk} }imply that
\begin{eqnarray}
f(x_k)-f(x_{k+1})
& \geq & [\barf(x_k) - \barf(x_{k+1})]
         - 2 \nu\barDt_{f,j_k}(x_k,s_k) \neol
& \geq & \eta\barDt_{f,j_k}(x_k,s_k)
         -2\nu\barDt_{f,j_k}(x_k,s_k)\neol
&   =  & (\eta-2\nu)\barDt_{f,j_k}(x_k,s_k)\neol
& \ge  & \frac 1 2  \eta \barDt_{f,j_k}(x_k,s_k){ \ge 0.}\label{DF-1F}
\end{eqnarray}
Thus
\beqn{myeq0}
\indicS\indicF = \indicS\indicF\indic{\{\DF_{k,+} \ge 0\}},
\eeqn
where $\DF_{k,+} =  f(X_k)-f(X_{k+1})$. Moreover, if $\calM_k$ also occurs
with $r_k\ge \bar r$ (i.e., if $\indicS\indicE\indicbL=1$) and
\req{approx-min} fails for $x_k$ and $\ord_k$, we may then use
\eqref{trqne-model-decrease2} to deduce from \req{DF-1F} that 
\beqn{fdec2}
f(x_k)-f(x_{k+1})
\ge \frac{\varsigma\eta}{2q!}\left(\kappa_\ord\epsmin\right)^{q+1} > 0,
\eeqn
which implies that, as long as \req{approx-min} fails,
\beqn{myeq1}
\indicS\indicE\indicbL
= \indicS\indicE\indicbL\indic{\{\DF_{k,+} > 0\}}.
\eeqn
\item[\textit{ii)}] If iteration $k$ is unsuccessful, the mechanism of the
  \algn algorithm guarantees that $x_k=x_{k+1}$ and, hence, that
  $f(x_{k+1})=f(x_k)$, giving that $(1-\indicS)\DF_{k,+} =0$.
\end{itemize}
Setting $f_0\eqdef f(X_0)$ and using this last relation and AS.2, we have that, for
any $\ell\in\{0,...,N_{\epsilon}-1\}$, 
\beqn{diff-f}
f_0-\flow
\ge f_0-f(X_{\ell+1})
= \sum_{k=0}^{\ell}\DF_{k,+}
= \sum_{k=0}^{\ell}\indicS \DF_{k,+}
\eeqn
Remembering that $X_0$ and thus $f_0$ are deterministic and
taking the total expectation on both sides of \req{diff-f} then gives that
{
\beqn{myeq6}
f_0-\flow
= \expect[f_0-\flow]
\ge \sum_{k=0}^{\ell} \expect\Big[\indicS \DF_{k,+}\Big]
= \sum_{k=0}^{\ell} \expect\Big[ \expect\Big[ \indicS \DF_{k,+} \mid \calA_{k-1} \Big]\Big].
\eeqn
}
Now, for $k \in \iiz{\ell}$,
\[
\indicS \DF_{k,+}
= \indicS \indicF \DF_{k,+} + \indicS (1-\indicF)\DF_{k,+}
\]
and so, using the second part of \req{prob-ass},
{
\begin{align}
\expect\Big[\indicS \DF_{k,+} \mid \calA_{k-1}\Big]
& = \expect\Big[ \indicS \indicF \DF_{k,+} \mid \calA_{k-1}\Big]
   + \expect\Big[ \indicS (1-\indicF) \DF_{k,+} \mid \calA_{k-1}\Big] \neol
& \ge  \expect\Big[ \indicS \indicF \DF_{k,+} \mid \calA_{k-1}\Big].
\label{myeq2}
\end{align}
}
Thus, again using the law of total expectations, \req{myeq6} yields that
{
\beqn{myeq7}
f_0-\flow
\ge \sum_{k=0}^{\ell} \expect\Big[ \expect\Big[ \indicS \indicF \DF_{k,+} \mid
    \calA_{k-1}\Big]\Big]
= \sum_{k=0}^{\ell} \expect\Big[ \indicS \indicF \DF_{k,+} \Big].
\eeqn
}
Moreover, successively using \req{myeq0}, \req{E-def}, \req{myeq1} and \req{fdec2},
\begin{align}
\indicS \indicF \DF_{k,+}
&= \indicS \indicF \indic{\{\DF_{k,+} > 0\}}\DF_{k,+}\neol
&= \indicS \indicF \indicM  \indic{\{\DF_{k,+} > 0\}}\DF_{k,+}
     + \indicS \indicF (1-\indicM) \indic{\{\DF_{k,+} > 0\}} \DF_{k,+}\neol
&\ge \indicS \indicE\indic{\{\DF_{k,+} > 0\}}\DF_{k,+}\neol
&\ge \indicS \indicE \indicbL\indic{\{\DF_{k,+} > 0\}}\DF_{k,+} \neol
&\ge \frac{\varsigma\eta}{2q!}\left(\kappa_\ord\epsmin\right)^{q+1}
     \Big(\indicS \indicE \indicbL\Big). \label{myeq4}
\end{align}
Substituting this bound in \req{myeq7} then gives that, as long as \req{approx-min} fails for
iterations $\ii{\ell}$,
\beqn{myeq3}
f_0-\flow
\ge \frac{\varsigma\eta}{2q!}\left(\kappa_\ord\epsmin\right)^{q+1}
      \sum_{k=0}^{\ell}\expect\Big[\indicS \indicE \indicbL\Big].
\eeqn
We now notice that, by Definition \ref{def2},
\[
N_{AS}-1 \le \sum_{k=0}^{N_{\epsilon}-2} \indicS \indicE \indicbL,
\]
and therefore
\beqn{myeq5}
\expect[N_{AS}-1]
\le \sum_{k=0}^{N_{\epsilon}-2}\expect\Big[\indicS \indicE \indicbL \Big].
\eeqn
Hence, letting $\ell=N_{\epsilon}-2$, substituting \req{myeq5} in \req{myeq3}, we deduce that  
\[
\expect[N_{AS}-1]\frac{\varsigma\eta}{2q!}\left(\kappa_\ord\epsmin\right)^{q+1}
\le f_0-\flow 
\]
and \req{boundEN1} follows.
} 

\noindent
While inequalities \req{boundEU} and \req{boundEN1} provide upper bounds on
$\expect[N_{AS}]$ and $\expect[N_U]$, as desired, the first still 
depends on $\expect[N_{IS}]$, which has to be bounded from above as well. This
can be done by following \cite{CartSche17} once more: 
Definition \ref{def2}, \eqref{EM1} and \eqref{defEM2} directly imply that 
\begin{equation}
\label{1boundEM3}
\expect[N_{IS}]
\leq \expect[N_I]
\leq  \frac{1-\pM}{\pM} \expect[N_A]
\leq \frac{1-\pM}{\pM}\left(\expect[N_{AS}] +  \expect[N_U] \right)
\end{equation}
and hence
\begin{equation}
\label{boundEM3}
\expect[N_{IS}]
\leq \frac{1-\pM}{2\pM-1}\left(2\expect[N_{AS}]
      +\left\lceil\log_{\gamma}\left( \frac{r_0}{\rr} \right)\right\rceil\right)
\end{equation}
follows from \eqref{boundEU} (remember that $\half < \pM \leq 1$). Thus, the
right-hand side in \eqref{EM1} is in turn bounded above because of \req{defEM2}, \req{boundEU},
\req{boundEM3} and \req{boundEN1}, giving
\begin{eqnarray}
\expect[N_A]
& \leq & \expect[N_{AS}]+\expect[N_U]
   \leq 2\expect[N_{AS}]+\expect[N_{IS}]+\left\lceil\log_{\gamma}\left( \frac{r_0}{\rr} \right)\right\rceil\neol
& \leq & \left(\frac{1-\pM}{2\pM-1}+1\right)\left(2\expect[N_{AS}]
      +\left\lceil\log_{\gamma}\left( \frac{r_0}{\rr} \right)\right\rceil\right)\neol
& \leq & \frac{\pM}{2\pM-1}\left[ \frac{4q!(f_0-\flow)}{ \varsigma\eta\left(\kappa_\ord\epsmin\right)^{q+1}} 
        +\left\lceil\log_{\gamma}\left( \frac{r_0}{\rr} \right)\right\rceil+2\right].\label{EM2}
\end{eqnarray}
This inequality, together with \eqref{PlanEN} and \eqref{EM1}, finally
gives the desired bound
\begin{equation}
\label{bound_Nsc}
\expect[N_{\Lambda}]
\leq \frac{1}{\pM}\expect[N_A]
\leq \frac{1}{2\pM-1}\left[ \frac{4q!(f_0-\flow)}{ \varsigma\eta\left(\kappa_\ord\epsmin\right)^{q+1}} 
        +\left\lceil\log_{\gamma}\left( \frac{r_0}{\rr} \right)\right\rceil+2\right].
\end{equation}

\noindent
We can now express our final complexity result in full.

\lthm{complexity-th}{
Suppose that AS.1--AS.3 hold, then
\beqn{trqne-complexity}
\expect[N_{\epsilon}]\le \frac{2\pM}{(2\pM-1)^2}\left[ \frac{4q!(f_0-\flow)}{ \varsigma\eta\left(\kappa_\ord\epsmin\right)^{q+1}} 
        +\left\lceil\log_{\gamma}\left( \frac{r_0}{\rr} \right)\right\rceil+2\right],
\eeqn
with $N_\epsilon$, $\rr$ and $\kappa_\ord$ defined as in
\req{Neps-def}, \eqref{trqne-condsucc} and  \eqref{trqne-kappddef2}, respectively.
}

\noindent
\proof{
Recalling the definitions \eqref{defnsigma} and the bound \eqref{bound_Ns}, we obtain that
\[
\expect[N_{\epsilon}]=\expect[N_{\Lambda}^c]+\expect[N_{\Lambda}]
\leq \frac{\expect[N_{\epsilon}]}{2\pM}+\expect[N_{\Lambda}],
\]
which implies, using \req{bound_Nsc}, that
\[
\frac{2\pM-1}{2\pM}\expect[N_{\epsilon}]
\le \frac{1}{2\pM-1}\left[ \frac{4q!(f_0-\flow)}{\varsigma\eta\left(\kappa_\ord\epsmin\right)^{q+1}} 
        +\left\lceil\log_{\gamma}\left( \frac{r_0}{\rr} \right)\right\rceil+2\right].
\]
This bound and the inequality $\half < \pM \leq 1$ yield the desired result.
} 

\noindent
We note that the $\calO\left(\epsmin^{-(q+1)}\right)$ evaluation bound given by
\req{trqne-complexity} is known to be sharp in order of $\epsmin$ for
trust-region methods using exact evaluations of functions and derivatives (see
\cite[Theorem~12.2.6]{CartGoulToin21}), which implies that
Theorem~\ref{complexity-th} is also sharp in order.

We conclude this section by noting that alternatives to the second
part of \req{prob-ass} do exist. For instance, we could assume that
\[
\expect\Big[\indicS \DF_{k,+} \mid \calA_{k-1}\Big]
\ge  \mu \expect\Big[ \indicS \indicF \DF_{k,+} \mid \calA_{k-1} \Big]
\]
for some $\mu>0$.  This condition can be used to replace the second part
of \req{prob-ass} to ensure \req{myeq2} in the proof of Lemma~\ref{LemmaNAS} and all
subsequent arguments.

\numsection{The impact of noise for first-order minimization}\label{noise-s}

While the above theory covers inexact evaluations of the objective function
and its derivatives, it does rely on AS.3.  Thus, as long as
the inexactness/noise on these values remains small enough for this assumption
to hold, the \al{TR$q$NE} algorithm iterates ultimately
produce an approximate local minimizer.  There are however practical
applications, such as minimization of finite sum using sampling strategies
(discussed in more detail below), where AS.3 may be unrealistic because of
noise intrinsic to the application.  We already
saw that, under the assumptions of Theorem~\ref{expect-ass-th}, a large enough
value of $\barDT_{f,j_k}(X_k,S_k)$ is sufficient for ensuring the third
condition in AS.3, but we also know from \req{Fk},\req{Ek}, \req{Df-vs-DT} and the definition of
$\nu$ that, at successful iterations,
\[
\Df(X_k,S_k)
\geq \eta \barDT_{f,j_k}(X_k,S_k)-E_k
\geq (\eta-2\nu)\barDT_{f,j_k}
\geq \half \eta \barDT_{f,j_k}
\]
whenever $\calF_k$ holds. Thus a large $\barDT_{f,j_k}(X_k,S_k)$ is only 
possible if either $\Df(X_k,S_k)$ is large or $\calF_k$ fails.  But a large
$\Df(X_k,S_k)$ is impossible close to a
(global) minimizer, and thus either $\calF_k$ (and AS.3) fails, or the guarantee
that the third condition of AS.3 holds vanishes when approaching a minimum.

Clearly, the above theory does not say anything about what happens for the
algorithm once AS.3 fails due to intrinsic noise.  Of course, this does not
mean it will not proceed meaningfully, but we can't guarantee it.
In order to improve our understanding of what can happen, we need to
consider particular realizations of the iterative process where AS.3 fails.
This is the objective of this section  where we focus on the instantiation 
\al{TR$1$NE} of \al{TR$q$NE} for first-order optimality.
Fortunately, limiting one's ambition to first order results in subtantial
simplifications in the \algn algorithm. 
We first note that the mechanism of Step~1 of \algn (whose purpose is to
determine $j_k$) is no longer necessary since $j_k$ must be equal to one if
only (approximate) gradients are available, so we can implicitly set
\[
D_{k,1}=\frac{\overline{\nabla_x^1 f}(X_k)}{\|\overline{\nabla_x^1 f}(X_k)\|}\Delta_k
\mand
\Delta T_{f,1}(X_k,D_{k,1})
= -\overline{\nabla_x^1 f}(X_k)^TD_{k,1}
=\|\overline{\nabla_x^1 f}(X_k)\|\Delta_k
\]
and immediately branch to the step computation. This in turn simplifies
to
\[
S_k=\frac{\overline{\nabla_x^1 f}(X_k)}{\|\overline{\nabla_x^1 f}(X_k)\|}R_k
\mand
\Delta T_{f,1}(X_k,S_k)
=\|\overline{\nabla_x^1 f}(X_k)\|\,R_k
\]
irrespective of the value of $\theta$, and \req{trqne-decrease2} automatically
holds. We thus observe that the simplified algorithm no longer needs
$\delta_{k,j}$ (since neither $\barphi_{f,j}^{\delta_k}(x_k)$ or
$\barDt_{f,j}(x_k,d_{k,j})$ needs to be effectively calculated), that the computed step $s_k$ is the global minimizer within the trust-region and that the
constant $\theta$ (used in Step~1 and the start of Step~2 of the \algn algorithm) is no
longer necessary. The resulting streamlined \al{TR$1$NE} algorithm is stated
as Algorithm~\ref{algo1} \vpageref{algo1}.

\algo{algo1}{The TR{\boldmath 1}NE algorithm}
{
\begin{description}
\item[Step~0: Initialisation.]
  A starting point $x_0$, a maximum radius $r_{\max}>0$ and an accuracy
  level $\epsilon \in (0,1)$ are given. The initial trust-region radius
  $r_0\in (\epsilon,  r_{\max}]$ is also given.
	For a given constant $\eta
  \in (0,1)$, define $  \nu \eqdef \min\big[\half \eta, \quarter(1-\eta)\big]$.
	Set $k=0$.

\vskip 5pt
\item[Step~1: Derivatives estimation.]
	Compute the derivative estimate $\overline{\nabla_x f}(x_k)$.  
	 \vskip 5pt
 \item[Step~2: Step computation.]
  Set $s_k =-\displaystyle \frac{\overline{\nabla^1_x f}(x_k)}{\|\overline{\nabla^1_x f}(x_k)\|} r_k$.
 \item[Step~3: Function decrease estimation.] Compute the estimate $\barf(x_k)-\barf(x_k+s_k)$ of
   $f(x_k)-f(x_k+s_k)$.
    \vskip 5pt
\item[Step~4: Test of acceptance.] Compute
  $ 
  \rho_k = \bigfrac{\barf(x_k) - \barf(x_k+s_k)}
                {\|\overline{\nabla^1_x f}(x_k)\|r_k}.
  $\\ 
  If $\rho_k \geq \eta$ \textit{(successful iteration)}, then set
  $x_{k+1} = x_k + s_k$; otherwise \textit{(unsuccessful iteration)} set $x_{k+1} = x_k$.
  \vskip 5pt
\item[Step~5: Trust-region radius update.]
  Set
  \[
  r_{k+1} = \left\{ \begin{array}{ll}
  {}\frac{1}{\gamma}r_k, & \tim{if} \rho_k < \eta\\
  {} \min[r_{\max},\gamma r_k], & \tim{if} \rho_k \ge  \eta.
  \end{array}\right.
  \]
  Increment $k$ by one and go to Step~1.
  \end{description}
  }

\noindent

  
The definition of the event $\calM_k$  in \req{Mk} ensures that
$\AM_k^{(2)}$ implies $\AM_k^{(1)}$ when first-order models are considered,
and thus, using also \eqref{Fk}, that $\calM_k$ and $\calF_k$ then reduce to 
\[
\AM_k =\{ | \|\nabla_x^1 f(X_k)\| -\|\overline{\nabla_x^1 f}(X_k)\| |
\le \nu  \|\overline{\nabla_x^1 f}(X_k)\| \}
\] 
and
\[
\calF_k =\{ |\Delta F(X_k,S_k) -{\Delta  f}(X_k,S_k)|
\le 2 \nu  \|\overline{\nabla_x^1 f}(X_k)\| R_k\},
\]
respectively.  Observe now that, because of the triangle inequality, $\calM_k$
is true whenever the event
\beqn{tildeM}
\widetilde{\calM}_k  \eqdef \{ \|\nabla_x^1 f(X_k)-\overline{\nabla_x^1 f}(X_k)\|
\le \nu  \|\overline{\nabla_x^1 f}(X_k)\| \}
\eeqn
holds, and, since
$\nu  \|\overline{\nabla_x^1 f}(X_k)\| \min \{1,R_k\}\le  \nu \|\overline{\nabla_x^1 f}(X_k)\|$, 
it also follows that $\calF_k$ is true whenever the event 
\begin{equation}\label{tildeF}
\widetilde{\calF}_k \eqdef \{| \Delta F(X_k,S_k) -{\Delta f}(X_k,S_k)|
\le 2 \nu  \|\overline{\nabla_x^1 f}(X_k)\| \min \{1,R_k\}\}
\end{equation}
holds. As a consequence,
\beqn{probF_tildeF}
\prob\Big [\calM_k \mid  \calA_{k-1}\Big ] \ge \prob\Big [\widetilde{\calM}_k \mid \calA_{k-1}\Big ]
\tim{and}
\prob\Big [\calF_k \mid  \calA_{k-1}\Big ] \ge \prob\Big [\widetilde{\calF}_k \mid \calA_{k-1}\Big ].
\eeqn
Our analysis of the impact of noise on the \al{TR$1$NE} algorithm starts by
considering a relatively general form for error distributions (as we did 
in Theorem~\ref{expect-ass-th}) and we then specialize our arguments to the
particular case of finite sum minimization with subsampling.

\subsection{Failure of AS.3 for general error distributions}

\noindent
At a generic iteration $k$,
suppose that $H_{0,k}$ and $H_{1,k}$, are continuous and increasing random functions
from $\Re_+$ to $[0,1]$ which are measurable for $\calA_{k-1}$ and such that $H_{0,k}(0)=H_{1,k}(0)=0$,
$\lim_{\tau \rightarrow +\infty} H_{0,k}(\tau)=\lim_{\tau \rightarrow +\infty}H_{1,k}(\tau)=1$ and,

\begin{eqnarray}
\prob\Big [| \Delta F(X_k,S_k) -{\Delta  f}(X_k,S_k)|<  \tau | \calA_{k-1}\Big] &\ge& H_{0,k}( \tau) \neol
\prob \Big [\|\nabla_x^1 f(X_k) -\overline{\nabla_x^1 f}(X_k)\|<  \tau | \calA_{k-1} \Big ] &\ge& H_{1,k}( \tau)\label{eq4.5}
\end{eqnarray}
For sake of simplicity, assume  $\alpha_*=\gamma_*\ge \sqrt{\half}$  in AS.3
and let $B_0$ and $B_1$  such that $H_{0,k}(B_0)=\sqrt{\alpha_*} $ and 
$H_{1,k}(B_1)= \sqrt{\alpha_*}$, and  $B=\max[B_0,B_1]$. Then,
\begin{align}
\prob\Big [| \Delta F(X_k,S_k) -{\Delta  f}(X_k,S_k)|<  \tau \mid \calA_{k-1}, \tau \ge B \Big ]
  \ge \sqrt{\alpha_*}, && \label{prb1}\\
\prob\Big [\|\nabla_x^1 f(X_k) -\overline{\nabla_x^1 f}(X_k)\|<  \tau \mid \calA_{k-1}, \tau \ge B  \Big]
\ge \sqrt{\alpha_*} . \label{prb2}
\end{align}
Define
\beqn{bar_beta}
\bar B=\frac{B}{\nu\min\{1,R_k\}}\ge \frac{B}{\nu} > B,
\eeqn
and note that $\bar B$ is measurable for $\calA_{k-1}$. Then \req{prb1} and \req{prb2} ensure that
 \begin{eqnarray}
   \prob\Big [| \Delta F(X_k,S_k) -{\Delta  f}(X_k,S_k)|\le \bar B
     \mid \calA_{k-1} \Big] \ge \sqrt{\alpha_*} \label{prob_f}\\
  \prob\Big[\|\nabla_x^1 f(X_k) -\overline{\nabla_x^1 f}(X_k)\|\le \bar B
    \mid \calA_{k-1} \Big] \ge \sqrt{\alpha_*}.  \label{prob_g}
  \end{eqnarray}
Finally, define the events
  \begin{align}
  \calG_k &\eqdef \{ \|\nabla_x^1 f (X_k)\|\ge  2\bar B \},  \label{eventG}\\
  \bar{\calG}_k &\eqdef \{\|\overline{\nabla_x^1 f}(X_k)\|\ge \bar B \}, \label{eventbarG}\\
  \calV_k &\eqdef \{ \|\nabla_x^1 f(X_k) -\overline{\nabla_x^1 f}(X_k)\|< \frac{B}{\nu}\} ,\label{eventV}
 \end{align} 
and observe that \req{prb2} implies that
 \beqn{probV}
 \prob\Big[\calV_k \mid \calA_{k-1}\Big ]\ge\sqrt{\alpha_*}.
 \eeqn
 
\lthm{barg-vs-g}{
 Let $\bar B$ as in \eqref{bar_beta}.  Then,  for each iteration $k$ of
 the \al{TR$1$NE} algorithm,
\beqn{boundprobbarG}
\prob\Big[\bar{\calG}_k \mid \calA_{k-1}, \calG_k \Big ] \ge \sqrt{\alpha_*}.
\eeqn
}

\proof{
For any realization of the \al{TR$1$NE} algorithm  we have that
$$
\|\overline{\nabla_x^1 f}(x_k)\|\ge \left| \|\nabla_x^1 f(x_k)  \| -\|\nabla_x^1 f(x_k) -\overline{\nabla_x^1 f}(x_k)\|   \right|.
$$
Therefore,  $ \|\nabla_x^1 f(x_k)\|\ge 2\bar \beta $ (where $\bar
\beta$ is the realization of $\bar B$) and
$\|\nabla_x^1 f(x_k) -\overline{\nabla_x^1 f}(x_k)\| \le  \beta/\nu <\bar\beta$ ensure that
$
\|\overline{\nabla_x^1 f}(x_k)\|\ge \bar \beta.
$ 
Then,  $\calG_k \cap \calV_k$ implies $\bar{\calG}_k$, where the events
$\calG_k$, $\calV_k$ and $\bar{\calG}_k$ are defined in
\eqref{eventG}-\eqref{eventV}, and  
$\prob\Big [\bar{\calG}_k \mid \calA_{k-1}, \calG_k, \calV_k \Big ] =1$. We
therefore have that
\[
\expect\Big[\indic{\bar{\calG}_k} \mid \calA_{k-1}, \calG_k \Big ]
\ge \expect\Big[\indic{\bar{\calG}_k} \mid \calA_{k-1}, \calG_k , \calV_k \Big ] \,
  \prob\Big[ \calV_k \mid \calA_{k-1}, \calG_k \Big]
  \ge   1 \cdot \sqrt{\alpha_*},
  \]
where we have used \eqref{probV} and the fact that
\[
\expect\Big[\indic{\calV_k} \mid \calA_{k-1}, \calG_k \Big]
= \expect\big[\expect\big[\indic{\calV_k} \mid \calA_{k-1} \big] \mid \calA_{k-1}, \calG_k \Big]
\]
to derive the last inequality. The conclusion \req{boundprobbarG} then follows.
} 

\lthm{mcondition_fail}{
Let $\bar B$ be defined by \req{bar_beta}.  Then,  for each iteration $k$ of
the \al{TR$1$NE} algorithm,
\begin{equation} \label{boundM}
\prob\Big [ \AM_k  \mid \calA_{k-1}, \left\{ \|\nabla_x^1 f(X_k)\|\ge 2 \bar B\right\}
\Big ] \ge  \alpha_*.
\end{equation}
Moreover, if $\omega$ is a realization for which
$\prob\Big [ \AM_k  \mid \calA_{k-1}\Big ]\rlz < \alpha_*$,  then
\beqn{degrad1}
\|\nabla_x^1 f(x_k)\| < 2\bar B\rlz.
\eeqn
}

\proof{
We obtain from \req{probF_tildeF} and  \req{boundprobbarG} that
\begin{eqnarray}
 \prob\Big[ \AM_k \mid \calA_{k-1}, \calG_k \Big] & \ge &\prob\Big[ \tAM_k \mid \calA_{k-1}, \calG_k \Big]
 =\expect\Big[\indictM \mid \calA_{k-1}, \calG_k \Big] \neol
&\ge& \expect\Big[ \indictM \mid \calA_{k-1}, \calG_k, \bar{\calG}_k \Big]\,
     \prob\Big[ \bar{\calG}_k \mid \calA_{k-1}, \calG_k \Big] \neol
& \ge&  \sqrt{\alpha_*} \,\expect\Big[\indictM \mid \calA_{k-1}, \calG_k, \bar{\calG}_k \Big].\label{boundMM}
\end{eqnarray}
If $\bar{\cal G}_k$ is true, then it follows from \eqref{bar_beta} and \req{eventbarG} that
$
\nu  \|\overline{\nabla_x^1 f}(X_k)\| \ge B.
$ 
Then, \eqref{tildeM} and \eqref{prb2} yield
\beqn{bounftildeM}
\expect\Big[ \indictM \mid \calA_{k-1}, \bar{\calG}_k \Big] \ge \sqrt{\alpha_*}.
\eeqn
Because the trace $\sigma$-algebra $\{ \calA_{k-1}, \bar{\calG}_k\}$ contains
the trace $\sigma$- algebra
$\{\calA_{k-1}, \calG_k, \bar{\calG}_k\}$, the tower property and
\eqref{bounftildeM} then imply that
\[
\begin{array}{lcl}
\expect\Big[\indictM \mid \calA_{k-1}, \calG_k, \bar{\calG}_k \Big]
& = &  \expect\Big[\expect\Big[ \indictM \mid \calA_{k-1}, \bar{\calG}_k \Big] \mid 
        \calA_{k-1}, \calG_k, \bar{\calG}_k \Big]\\*[2ex]
& \ge & \expect\Big[ \sqrt{\alpha_*} \, \mid \calA_{k-1}, \calG_k, \bar{\calG}_k \Big]
 =  \sqrt{\alpha_*}
\end{array}
\]
which, together with \req{boundMM} gives \eqref{boundM}.
Since $\calG_k$ is measurable for $\calA_{k-1}$ we have 
\[
\prob\Big[\AM_k \mid \calA_{k-1}\Big]
\ge \prob\Big[\AM_k \mid \calA_{k-1},\calG_k \Big] \expect\Big[\indic{\calG_k} \mid \calA_{k-1}\Big]
 = \prob\Big[ \AM_k \mid \calA_{k-1},\calG_k\Big]\, \indic{\calG_k}.
 \]
If we now consider a realization $\omega$ such that
$\prob\Big[\AM_k \mid  \calA_{k-1}\Big]\rlz < \alpha_*$,
we therefore obtain, using \req{boundM} taken for the realization $\omega$, that
\[
\alpha_*
> \prob\Big[ \AM_k \mid \calA_{k-1},\calG_k\Big]\rlz \, \indic{\calG_k\rlz}
\ge \alpha_* \, \indic{\calG_k\rlz},
\]
which implies that $\indic{\calG_k\rlz}=0$, and thus that \req{degrad1} holds.
}

\lthm{fcondition_fail}{
 Let $\bar B$ be defined by \req{bar_beta}.  Then, for each iteration $k$ of
 the \al{TR$1$NE} algorithm,
\begin{equation} \label{boundF}
\prob\Big[ \calF_k  \mid \calA_{k-1},
\left\{ \|\nabla_x^1 f(X_k)\| \ge 2\bar B\right\} \Big]
\ge \alpha_*.
\end{equation}
Moreover, if $\omega$ is a realization for which
$\prob\Big[ \calF_k \mid \calA_{k-1}\Big ]\rlz  < \alpha_*$, then
\beqn{degrad2}
  \|\nabla_x^1 f(x_k)\|< 2\bar B\rlz.
\eeqn
}

\proof{The proof is similar to that of Theorem~\ref{mcondition_fail}, and is
  given in appendix for completeness.}

\noindent
Theorems~\ref{mcondition_fail} and \ref{fcondition_fail} indicate that the
assumptions made in AS.3 about $\calM_k$ and $\calF_k$ are likely to be
satisfied as long as the gradients remain sufficiently large, allowing the
\al{TR$1$NE} algorithm to iterate meaningfully.  Conversely,
they show that, should these assumptions fail for a particular realization of
the algorithm because of a high level of intrinsic noise, ``degraded''
versions of first-order optimality conditions given by \req{degrad1} and
\req{degrad2} nevertheless hold when this failure occurs.

\subsection{A subsampling example}

We finally illustrate how intrinsic noise might affect our probabilistic
framework on an example.
Suppose that
\beqn{finite-sum}
f(x) = \frac{1}{m}\sum_{i=1}^m f_i(x),
\eeqn
where the $f_i$ are functions from $\Re$ to $\Re$ having Lipschitz continuous
gradients and where $m$ is so large that computing the complete value of
$f(x)$ or its derivatives is impractical.  Such a situation occurs for
instance in algorithms for deep-learning, an application of growing
importance. A well-known strategy to obtain approximations of the desired
values at an iterate $x_k$ is to sample the $f_i(x_k)$ and compute the sample
averages, that is 
\beqn{subsampled}
\barf(x_k) = \frac{1}{|\mathfrak{b}_0(x_k)|}\sum_{i\in \mathfrak{b}_0(x_k)}f_i(x_k),
\ms
\overline{\nabla_x^1 f}(x_k)
= \frac{1}{|\mathfrak{b}_1(x_k)|}\sum_{i\in \mathfrak{b}_1(x_k)}\nabla_x^1f_i(x_k),
\eeqn
where $\mathfrak{b}_0(x)$ and $\mathfrak{b}_1(x)$ are realizations of random
``batches'', that is randomly selected\footnote{With uniform distribution.} subsets of
$\ii{m}$. Observe that Step~3 of the \al{TR$1$NE} algorithm
computes the estimate $\barf(x_k)-\barf(x_k+s_k)$, which we assume, in the
context of \req{finite-sum}, to be
\beqn{subsampled2}
\barf(x_k) - \barf(x_k+s_k) = \frac{1}{|\mathfrak{b}_0(x_k)|}\sum_{i\in \mathfrak{b}_0(x)}(f_i(x_k)-f_i(x_k+s_k)),
\eeqn
(using a single batch for both the function estimates).  Observe that our choice to make $\mathfrak{b}_0$ and
$\mathfrak{b}_1$ dependent on $x_k$ implies that their random counterparts $\mathfrak{B}_0(X_k)$ and
$\mathfrak{B}_1(X_k)$ are measurable for $\calA_{k-1}$ (clearly we could have chosen a
more complicated dependence on the past of the random process). The mean-value
theorem then yields that, for some
$\{y_i\}_{i\in \mathfrak{b}_0(x_k)}$ in the segment $[x_k,x_k+s_k]$,
\[
|\barf(x_k) - \barf(x_k+s_k)|
\le \left(\frac{1}{|\mathfrak{b}_0(x_k)|}\sum_{i\in \mathfrak{b}_0(x_k)}\nabla_x^1 f_i(y_i)\right)\|s_k\|
\le r_k \max_{\stackrel{{\scriptstyle y\in[x_k,x_k+s_k]}}{\scriptstyle i\in \mathfrak{b}_0(x_k)}}\|\nabla_x^1 f_i(y)\|
\]

Note that one expects the right-hand side of this inequality
to be quite small when the trust-region radius is small or when convergence to
a local minimizer occurs and $\overline{\nabla_x^1 f}(x_k)$ is a reasonable
approximation of $\nabla_x^1 f(x_k)$. To simplify our illustration, we assume,
for the rest of this section, that there exists a constant $\kappa_f$
 such that for any $y \in \Re^n$, for every realization
$\mathfrak{b}_0(x_k)$, 
\[
r_k \max_{i\in \{1,\ldots,n\}}\|\nabla_x^1 f_i(y)\|
\le \kappa_f.
\]
Returning to the random process and using the Bernstein concentration
inequality, it  results from \cite[Relation (7.8)]{BellGuriMoriToin19} that,
for any $k$ and deterministic $\tau>0$,
\beqn{concentration}
\prob\Big[\DF(X_k,S_k)-\Df(X_k,S_k)>  \tau \Big ] \le e^{-W_0(\tau)}
\tim{where}
W_0(\tau) = \frac{\tau^2|\mathfrak{B}_0(X_k)|}{4\kappa_f(2\kappa_f+\sfrac{1}{3}\tau)}.
\eeqn
Similarly,
\beqn{concentration2}
\prob\Big[\| \nabla_x^1 F(X_k,S_k)-\nabla_x^1f(X_k,S_k)\| >  \tau \Big ] \le
\min\left[1,(n+1)e^{-W_1(\tau)}\right],
\ms
W_1(\tau) = \frac{\tau^2|\mathfrak{B}_1(X_k)|}{4\kappa_g(2\kappa_g+\sfrac{1}{3}\tau)},
\eeqn
for some constant $\kappa_g>0$. One also checks that, since $\mathfrak{B}_0(X_k)$ and
$\mathfrak{B}_1(X_k)$ are measurable for $\calA_{k-1}$, so are $W_0$ and $W_1$.
One then easily verifies that $W_0(\tau)$ is an increasing function of $\tau$,
and hence $e^{-W_0(\tau)}$ is  decreasing. Letting
$G_k(\tau)=1-e^{-W_0(\tau)}$, we immediately obtain that conditions \req{g0} and \req{g1}
hold. Let us now analyze condition \req{g2} and consider any realization
$\omega$, where $w_0(\tau) \eqdef W_0\rlz(\tau)$. Note that
$w_0(\tau)\ge |\mathfrak{b}_0(x_k)|^{\half} \tau$ when 
\beqn{phistar-def}
\tau
\ge \tau_*
\eqdef \frac{8\kappa_f^2}{|\mathfrak{b}_0(x_k)|^{\half}-\frac{4}{3} \kappa_f}
\eeqn
and $\frac{|\mathfrak{b}_0(x_k)|^{\half}}{\kappa_f}$ is large enough so that 
\beqn{B0size}
\frac{|\mathfrak{b}_0(x_k)|^{\half}}{\kappa_f}>\frac{4}{3}.
\eeqn
Hence $e^{-w_0(\tau)}\le e^{-|\mathfrak{b}_0(x_k)|^{\half} \tau}$ for all $\tau\ge \tau_*$
and 
\begin{eqnarray*}
  \int_0^\infty e^{-w_0(\tau)}\,d\tau
  \le \int_0^{\tau_*} e^{-w_0(\tau)}\, d\tau+\int_{\tau_*}^\infty e^{-|\mathfrak{b}_0(x_k)|^{\half} \tau}\, d\tau.
\end{eqnarray*}
In addition, since $e^{-w_0(\tau)}$ is decreasing and non-negative, we have that
\[
\int_0^{\tau_*} e^{-w_0(\tau)}\, d\tau \le \tau_* e^{-w_0(0)}= \tau_*.
\]
This bound and \req{phistar-def} then imply that
\beqn{badbound}
\int_0^\infty (1-g_k(\tau))\, d\tau
=\int_0^\infty e^{-w_0(\tau)}\, d\tau
\le  \tau_*+\frac{e^{-|\mathfrak{b}_0(x_k)|^{\half} \tau_*}}{|\mathfrak{b}_0(x_k)|^{\half}}
< +\infty,
\eeqn
proving that \req{g2} also holds for the arbitrary realization $\omega$.  We
may therefore apply Theorem~\ref{expect-ass-th} provided $\mathfrak{B}_0(X_k)$ is sufficiently
large\footnote{%
While the bound given by \req{badbound} is adequate for our proof, this
inequality can be pessimistic.  For instance, if we set
$\kappa_f=1$ and $|\mathfrak{b}_0(X_k)| = 2056$, the numerically computed
value of the left-hand side is 0.0556 while that of the right-hand side is 0.1818.
}
to ensure \req{B0size} and \req{phistar-def}, and conclude that, under these conditions,
\req{expectation} holds whenever
\[
\Dt_{f,1}(x_k,s_k)
\ge \frac{1}{\eta}\int_0^\infty (1-g_k(\tau)) \,d\tau.
\]
We can also apply the analysis in Section 4.1  with
\[
H_{0,k}(\tau) = 1 - e^{-W_0(\tau)}
\mand
H_{1,k}(\tau) = \max\left[ 0,1 - (n+1)e^{-W_1(\tau)}\right].
\]
A short calculation shows that
$B_0 = \calO(\kappa_f|\mathfrak{B}_0|^{-\half}|)$ and
$B_1 = \calO(\kappa_g|\mathfrak{B}_1|^{-\half}|)$, where $B_0$ and $B_1$ are
defined below \req{eq4.5}. Then, 
Theorems~\ref{mcondition_fail} and \ref{fcondition_fail}  hold with $\bar B=\calO\big (\frac{\max \{\kappa_f|\mathfrak{B}_0|^{-\half}|,\kappa_g|\mathfrak{B}_1|^{-\half}|\}}{\nu\max\{1,R_k\}}\big)$.

We finally illustrate the impact of intrinsic noise on the (admittedly ad-hoc)
problem of minimizing
\beqn{ex-def}
f(x) = \frac{1}{2m}\sum_{i\in \calZ}\left[\half x^2 + \half \alpha \,\sgn(i)\,e^{- x^2}\right]
\eqdef \frac{1}{2m}\sum_{i\in \calZ} f_i(x),
\eeqn
where $\alpha >0$ is a noise level and where
$\calZ = \iibe{-m}{m}\setminus\{0\}$ for some large integer $m$.  Suppose furthermore that the $\{f_i(x)\}_{i\in\calZ}$ and 
$\{\nabla_x^1 f_i(x)\}_{i\in\calZ}$ are computed by black-box routines, therefore
hiding their relationships.
Consider an iterate $x_k$ at the start of iteration $k$ of an arbitrary
realization of the \al{TR$1$NE} algorithm\footnote{With given $\nu$ and $\varsigma=1$.} applied
to this problem.  We verify that, for $i\in\calZ$,
\[
\nabla_x^1f_i(x_k) = x_k\left(1-\alpha \,\sgn(i)\,e^{- x_k^2}\right)
\]
and thus $\nabla_x^1 f(x_k) =x_k$ and, in view of
\req{phi-def},
\beqn{ex-phi}
\phi_{f,1}^{\ord_k}(x_k) = |x_k|\,\ord_k
\eeqn
for all $x_k$.  As a consequence, $x=0$ is the unique global minimizer of
$f(x)$.
Suppose, for the rest of this section, that
$\mathfrak{B}_{0,k} \eqdef \mathfrak{B}_0(x_k)= \mathfrak{B}_0(x_k+s_k)$,
$\mathfrak{B}_{1,k} \eqdef \mathfrak{B}_1(x_k)$, and that
$n_{0,k}^\mathfrak{B}$ and $n_{1,k}^\mathfrak{B}$, the cardinalities of
these two sets are known parameters. We deduce from \req{subsampled} that
\beqn{app-grad}
\overline{\nabla_x^1 f}(x_k)
= x_k \left(1- \alpha \Psi(\mathfrak{B}_{1,k}) e^{-x_k^2}\right)
\tim{ where }
\Psi(\mathfrak{B}) \eqdef \frac{1}{|\mathfrak{B}|}\sum_{i\in\mathfrak{B}}\sgn(i).
\eeqn
Thus $\Psi(\mathfrak{B})$ is a zero-mean random variable with values in
$[-1,1]$, depending on the randomly chosen batch
$\mathfrak{B}\subseteq\calZ$ of size $|\mathfrak{B}|$. Using the hypergeometric distribution, it is
possible to show that $|\Psi(\mathfrak{B})|$ is (in probability) a
decreasing function of $|\mathfrak{B}|$.

Moreover, the use of standard tail bounds \cite{Chva79} reveals that, for any
$t\in(0,1)$, 
\beqn{ex-distrib}
\prob\Big[|\Psi(\mathfrak{B}_{1,k})|\leq t \Big]
= \prob\Big[\Psi(\mathfrak{B}_{1,k})\leq t \Big]^2
= \left(1-\prob\Big[\Psi(\mathfrak{B}_{1,k}) > t \Big]\right)^2
\ge (1-e^{-\half t^2 n_{1,k}^\mathfrak{B}})^2,
\eeqn
in turn indicating that $\prob\Big[|\Psi(\mathfrak{B}_{1,k})| \le t \Big] >
\half$ whenever
\[
n_{1,k}^\mathfrak{B}
\ge \frac{2}{t^2}\left|\log\left(1-\frac{1}{\sqrt{2}}\right)\right|
\approx \frac{2.4559}{t^2}.
\]

\vspace*{2mm}
\noindent
\textbf{Occurence of $\calM_k^{(1)}$ and  $\calM_k^{(2)}$. }
Let us now examine at what conditions the events $\calM_k^{(1)}$ and
$\calM_k^{(2)}$ do occur for a specific realization $\mathfrak{b}_{1,k}$ of
$\mathfrak{B}_{1,k}$ , and consider the occurence of $\calM_k^{(1)}$ first.
Because the minimum of first-order models in a ball of radius $\ord_k$ must
occur on the boundary, we choose $d_{k,1} = -\sgn(\overline{\nabla_x^1f}(x_k))\ord_k$
so that
\[
\barDt_{f,1}(x_k,d_{k,1}) = |\overline{\nabla_x^1f}(x_k)| \, \ord_k.
\]
Using \req{app-grad}, we then have that
\beqn{Dtk}
\barDt_{f,1}(x_k,d_{k,1})
=|x_k|\,\ord_k\,\left|1 - \alpha\Psi(\mathfrak{b}_{1,k}) \,e^{- x_k^2}\right|.
\eeqn
Thus the quantity $1-\alpha\Psi(\mathfrak{b}_{1,k})\,e^{- x_k^2}$ may be
interpreted as the local noise relative to the model decrease.

\begin{figure}[htb]
\centerline{   
\includegraphics[width=6cm]{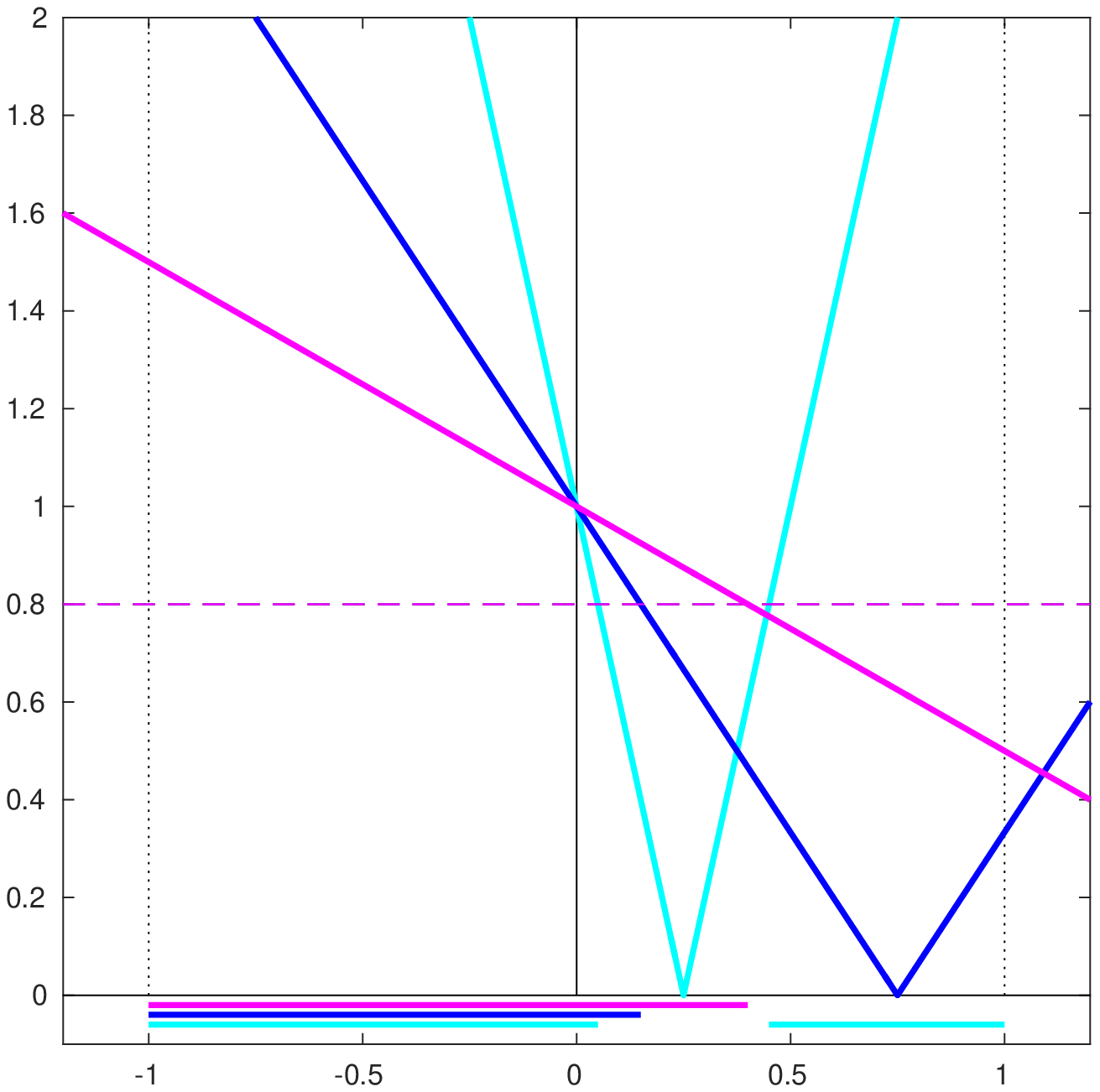}
\hspace*{10mm}
\includegraphics[width=6cm]{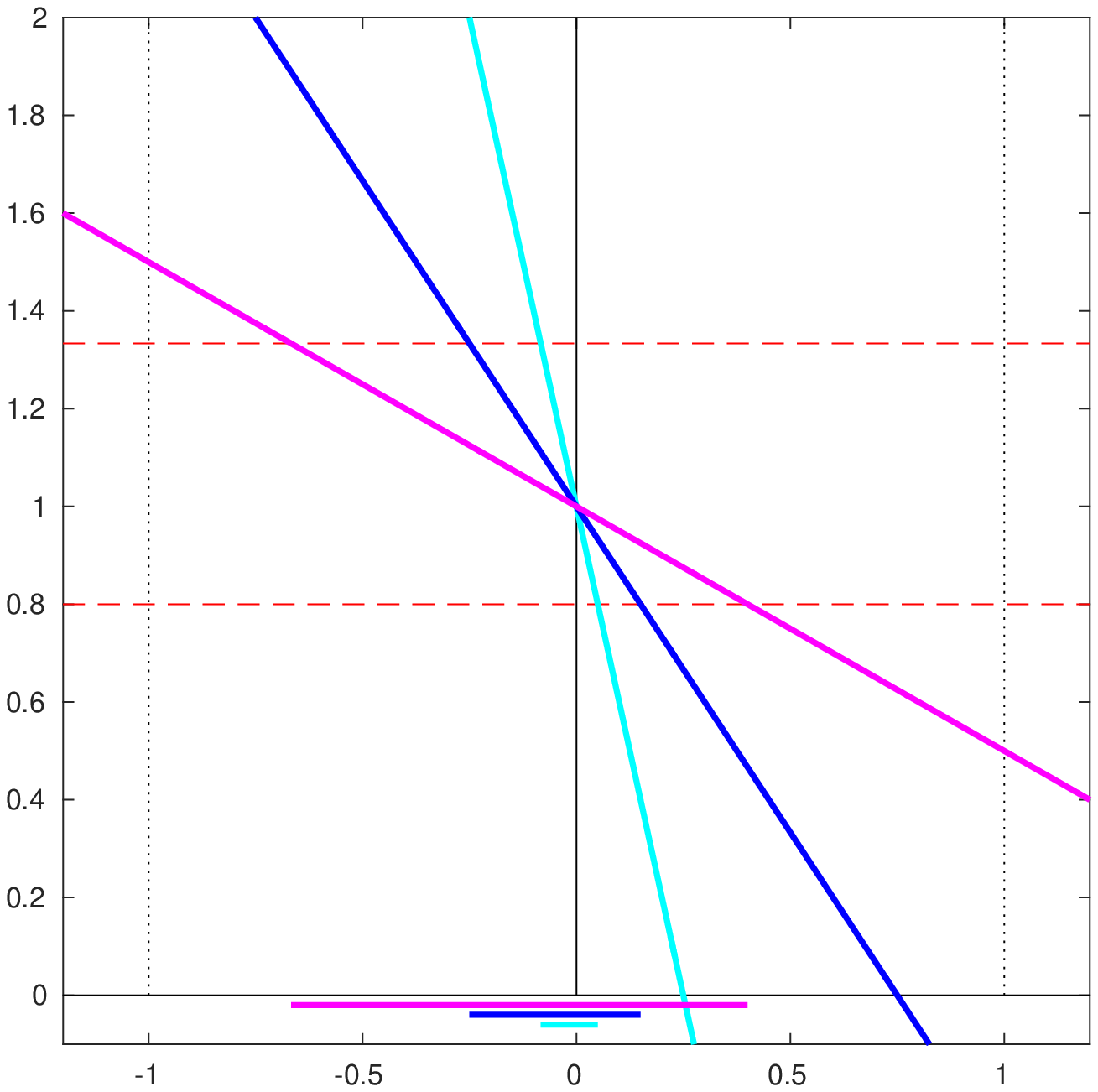}
}
\caption{\label{exfig}
An illustration of conditions \req{ex-condM1} (left) and
\req{ex-condM2} (right) as a function of $\Psi(\mathfrak{b}_{1,k})$,
for $\nu = \quarter$ and local relative noise levels $\alpha e^{-x_k^2}=
\half$ (magenta), $\sfrac{4}{3}$ (blue) and $4$ (cyan). Acceptable ranges for
$\Psi(\mathfrak{b}_{1,k})$ are shown below the horizontal axis in matching colours.
}
\end{figure}

Taking \req{ex-phi} and \req{Dtk} into account, $\calM_k^{(1)}$ occurs, in any
realization, whenever
\beqn{ex-condM1}
\left|1 -\alpha\Psi(\mathfrak{b}_{1,k})\,e^{- x_k^2}\right|
\ge \frac{1}{1+\nu},
\eeqn
that is
\[
\Psi(\mathfrak{b}_{1,k}) \le
\frac{e^{x_k^2}}{\alpha}\left(1-\frac{1}{1+\nu}\right)
\tim{ or }
\Psi (\mathfrak{b}_{1,k}) \ge
\frac{e^{x_k^2}}{\alpha}\left(1+\frac{1}{1+\nu}\right).
\]
This condition may be quite weak, as shown in left picture in
Figure~\ref{exfig}, where the shape of the left-hand side of \req{ex-condM1}
is shown for increasing values\footnote{Magenta for 0.5, blue for 4/3 and cyan
for 4.} of the local noise level $\alpha \,{\rm exp}(-x^2)$ as a function
of $\Psi(\mathfrak{b}_{1,k})$, and where the lower bound $1/(1+\nu)$ is shown as a
red horizontal dashed line. The corresponding ranges of acceptable values of
$\Psi(\mathfrak{b}_{1,k})$ are shown below the horizontal axis (in matching
colours). The one-sided nature of the inequality defining {$\calM_k^{(1)}$} is
apparent in the picture, where restrictions on the acceptable values of
$\Psi(\mathfrak{b}_{1,k})$ only occur for positive values. This reflects the fact that
the model may be quite inaccurate and yet produce a decrease which is large
enough for the condition to hold. 

The constraints on $\Psi(\mathfrak{b}_{1,k})$, and thus on $n_{1,k}^\mathfrak{B}$, become
more stringent when considering the occurence of $\calM_k^{(2)}$.  Since, for
any realization, 
$s_k = -\sgn\big(\overline{\nabla_x^1 f}(x_k)\big)r_k$,
we deduce from \req{app-grad} that
\begin{align*}
\Dt_{f,1}(x_k,s_k)
= -x_ks_k
&= -x_k \left[-\sgn(x_k)\,
  \sgn\left(1 - \alpha \,\Psi(\mathfrak{b}_{1,k})\,e^{- x_k^2}\right)r_k\right]\neol
& = |x_k| \, r_k \,\sgn\left(1- \alpha \,\Psi(\mathfrak{b}_{1,k})\,e^{- x_k^2}\right)
\end{align*}
and
\beqn{barDts}
\barDt_{f,1}(x_k,s_k)
=|x_k|\,r_k\,
\left|1-\alpha\Psi(\mathfrak{b}_{1,k})\,e^{- x_k^2}\right|.
\eeqn
One then verifies that $\calM_k^{(2)}$ occurs whenever
\beqn{ex-condM2}
\frac{1}{1+\nu}
\le 1-\alpha\Psi(\mathfrak{b}_{1,k})\,e^{- x_k^2}
\le \frac{1}{1-\nu}.
\eeqn
The acceptable values for $\Psi(\mathfrak{b}_{1,k})$ are illustrated in the right
picture of Figure~\ref{exfig}, which
{%
shows the shape of the central term in \req{ex-condM2} using
}%
the same conventions than for the
left picture except that now the acceptable part of the curves lies between
the lower and upper bounds resulting from \req{ex-condM2} (again shown as
dashed red lines).
{
A short calculation reveals that \req{ex-condM2} is equivalent to requiring
\[
\frac{e^{x_k^2}}{\alpha}\left(1-\frac{1}{1-\nu}\right)
\leq \Psi(\mathfrak{b}_{1,k})
\leq \frac{e^{x_k^2}}{\alpha}\left(1-\frac{1}{1+\nu}\right).
\]
}%
This therefore defines intervals around the origin, whose
widths clearly decrease with the local relative noise level.
Because $|\Psi(\mathfrak{B}_{1,k})|$ is (in probability) a decreasing function of
$n_{1,k}^\mathfrak{B}$, this indicates that $n_{1,k}^\mathfrak{B}$ must increase with
$\alpha e^{-x_k^2}$, that is when the local relative noise is large.

\vspace*{2mm}
\noindent
\textbf{Occurence of $\calF_k$. }
A similar reasoning holds when considering the event $\calF_k$.
Given \req{subsampled2}, we have that
\beqn{errDfex}
|\Df(x_k,s_k) - \barDf(x_k,s_k)|
= \half \alpha \,|\Psi(\mathfrak{b}_{0,k})|\,\left|e^{- x_k^2} - e^{-(x_k+s_k)^2}\right|
\eeqn
and, in view of \req{barDts}, $\calF_k$ thus occurs whenever
\beqn{ex-condF}
\half \alpha \,|\Psi(\mathfrak{b}_{0,k})|\,e^{-x_k^2}\left|1 - e^{-(2x_ks_k+s_k^2)}\right|
\le 2 \nu |x_k|\,r_k\,\left|1-\alpha\Psi(\mathfrak{b}_{1,k})\,e^{-x_k^2}\right|.
\eeqn
Thus, if $|x_k|$ is small ({e.g.}, if the optimum is close)
then satisfying \req{ex-condF} requires the left-hand side of this inequality
to be small, putting a high request on $n_{0,k}^\mathfrak{B}$, while the inequality
is more easily satisfied if $|x_k|$ is large, irrespective of the batch sizes.
Note that, in the first case (i.e., when $|x_k|$ is small), the request on
$n_{0,k}^\mathfrak{B}$ is stronger for smaller $n_{1,k}^\mathfrak{B}$.

\vspace*{2mm}
\noindent
\textbf{Occurence of \req{expectation}. }
Given \req{ex-distrib} and \req{barDts}, we see from 
Theorem~\ref{expect-ass-th} that \req{expectation} holds whenever
\[
\barDt_{f,1}(x_k,s_k)
=|x_k|\,r_k\,
\left|1-\alpha\Psi(\mathfrak{b}_{1,k})\,e^{- x_k^2}\right|
\ge \frac{1}{\eta} \int_0^\infty(1-g_k(\tau))\,d\tau
= \frac{1}{\eta}\sqrt{\frac{\pi}{2n_{1,k}^{\mathfrak{B}}}}
\]
where we have used the definition of the erf function to derive the last
equality. Thus, as $|\Psi(\mathfrak{b}_{1,k})|\le 1$, gauranteeing \req{expectation}
requires a larger $n_{1,k}^{\mathfrak{B}}$ for small value of $x_k$, that is
when the optimum is approached.

\numsection{Conclusions and perspectives}\label{section:conclusion}

We have considered a trust-region method for unconstrained minimization
inspired by \cite{CartGoulToin20c} which is adapted to handle randomly perturbed
function and derivatives values and is capable of finding approximate
minimizers of arbitrary order. Exploiting ideas of
\cite{CartSche17,BlanCartMeniSche19}, we have shown that its evaluation
complexity is (in expectation) of the same order in the requested accuracy as
that known for related deterministic methods
\cite{BlanCartMeniSche19,CartGoulToin20c}.

In \cite{BellGuriMoriToin21b}, the authors have considered the effect of
intrinsic noise on complexity of a deterministic, noise tolerant variant of
the trust-region algorithm.  This important question is handled here by considering
specific realizations of the algorithm under reasonable assumptions on the
cumulative distribution of errors in the evaluations of the objective function
and its derivatives. We have shown that, for such realizations, a first-order
version of our trust-region algorithms still provides ``degraded'' optimality
guarantees, should intrinsic noise cause the assumptions used for the
complexity analysis to fail. We have specialized and illustrated those results
in the case of sampling-based finite-sum minimization, a context of particular
interest in deep-learning applications.

We have so far developed and analyzed ``noise-aware'' deterministic and
stochastic algorithms for unconstrained optimization. Clearly, considering the
constrained case is a natural extension of the type of analysis presented
here.


\appendix

\appnumsection{Appendix: additional proofs}

\vskip 10 pt
\noindent
{\bf Proof of \req{CSlemma21}}
\vskip 5 pt
\noindent
\proof{Since  $\sigma(\indicLc)$
  belong to $\sigax$, because the random variable $\Lambda_k$ is fully determined, assuming $Pr(\indicLc)>0$,
the tower property yields:
$$  
\expect \left[ \indicE\mid \indicLc\right]=\expect\left [ \expect  \left[ \indicE \mid \calA_{k-1}\right]\mid \indicLc\right ]\ge \expect\left [ \pM \mid \indicLc\right ]=\pM.
$$
Then, by the total expectation law  we have
 $$
\expect \left[ \indicE \indicLc\right]=\expect\left [  \indicLc\expect\left[ \indicE \mid \indicLc\right ]\right ]\ge \pM \expect\left [ \indicLc \right ].
$$
Similarly,
 $$
\expect \left[ \indic{\{k<N_\epsilon\}}\indicE \indicLc\right]\ge \pM \expect\left [\indic{\{k<N_\epsilon\}} \indicLc \right ],
$$
as $\indic{\{k<N_\epsilon\}} $ is also determined by $\calA_{k-1}.$ In case $Pr(\indicLc)=0$,  the above inequality  holds trivially.
Then
$$
\expect\left[ \sum_{k=0}^{N_\epsilon-1}\indicLc\indicE \right]=\expect\left[ \sum_{k=0}^{\infty}\indic{\{k<N_\epsilon\}} \indicLc\indicE \right]\ge p_M
\expect\left[ \sum_{k=0}^{\infty}\indic{\{k<N_\epsilon\}} \indicLc \right ]=
\pM\,\expect\left[ \sum_{k=0}^{N_\epsilon-1}  \indicLc \right].
$$
and \req{CSlemma21} follows.}

\vskip 5 pt
\noindent
{\bf Proof of Lemma~\ref{a1}}
\vskip 5 pt
\noindent
\proof{Proceeding as in the proof of \req{CSlemma21} with $\indicL$ in place of $\indicLc$, we obtain:
$$
\expect\left[ \sum_{k=0}^{N_\epsilon-1}\indicL\indicE \right]\ge 
\pM\,\expect\left[ \sum_{k=0}^{N_\epsilon-1}  \indicL \right].
$$
Moreover, proceeding again as in the proof of \req{CSlemma21} and  substituting $\indicE$ with $\indicEc$ we obtain 
$$
\expect\left[ \sum_{k=0}^{N_\epsilon-1}\indicLc\indicEc \right]\le 
(1-\pM)\,\expect\left[ \sum_{k=0}^{N_\epsilon-1}  \indicL \right].
$$
Using the above inequalities we obtain \req{EM1}.}

\vskip 10 pt
\noindent
{\bf Proof of Theorem~\ref{fcondition_fail}}
\vskip 5 pt
\noindent
\proof{Because of \req{probF_tildeF} and \req{boundprobbarG}, we have that
\begin{eqnarray}
 \prob\Big[ \calF_k \mid  \calA_{k-1}, \calG_k \Big ]&\ge &\prob\Big[\widetilde{\calF}_k \mid \calA_{k-1}, \calG_k \Big ]
= \expect\Big [\indictF \mid \calA_{k-1}, \calG_k \Big ] \neol
&\ge& \expect\Big[\indictF \mid \calA_{k-1}, \calG_k, \bar{\calG}_k \Big] \,
    \prob\Big [\bar{\calG}_k\mid \calA_{k-1}, \calG_k \Big] \neol
&\ge& \sqrt{\alpha_*} \, \expect\Big[\indictF| \calA_{k-1}, \calG_k, \bar{\calG}_k \Big ].\label{boundFF}
\end{eqnarray}
If $\bar{\calG}_k$ is true then by \eqref{bar_beta} it follows
$$
2 \nu  \|\overline{\nabla_x^1 f}(X_k)\| \min \{1,R_k\}\ge \bar B \nu  \min \{1,R_k\}>B .
$$
Then, \eqref{tildeF} and \eqref{prb1} yield
\beqn{bounftildeF}
\expect\Big[\indictF \mid \calA_{k-\frac{1}{2}}, \bar{\cal G}_k\Big ] \ge \sqrt{\alpha_*}.
\eeqn
Because the trace $\sigma$-algebra $\{ \calA_{k-1}, \bar{\calG}_k\}$ contains
the trace $\sigma$- algebra
$\{\calA_{k-1}, \calG_k, \bar{\calG}_k\}$, the tower property and
\eqref{bounftildeF} then imply that
\[
\expect\Big[\indictF \mid \calA_{k-1}, \calG_k, \bar{\calG}_k \Big]
=  \expect\Big[\expect\Big[ \indictF \mid \calA_{k-1},\bar{\calG}_k \Big ] \mid
     \calA_{k-1}, \calG_k, \bar{\calG}_k\Big ]
\ge \expect\Big[\sqrt{\alpha_*} \, \mid \calA_{k-1}, \calG_k, \bar{\calG}_k \Big ] 
= \sqrt{\alpha_*}
\]
which, together with \req{boundFF}, implies \req{boundF}.
Since $\calG_k$ is measurable for $\calA_{k-1}$ we have that
\begin{eqnarray*}
  \prob\Big[ \calF_k \mid \calA_{k-1}\Big] \ge \prob\Big[\calF_k \mid \calA_{k-1},\calG_k \Big] \,\expect\Big[\indic{\calG_k} \mid \calA_{k-1}\Big]
= \prob\Big[ \calF_k \mid \calA_{k-1},\calG_k \Big] \indic{\calG_k}.\\
\end{eqnarray*}
Considering now a realization $\omega$ such that $\prob\Big[\calF_k \mid
  \calA_{k-1}\Big]\rlz <\alpha_*$, we therefore obtain, using \eqref{boundF} taken for this realization, that
\[
\alpha_*
> \prob\Big[ \calF_k \mid \calA_{k-1},\calG_k\Big]\rlz \,\indic{\calG_k\rlz}
\ge \alpha_* \, \indic{\calG_k\rlz},
\]
which implies that $\indic{\calG_k\rlz}=0$, in turn yielding \req{degrad2}.}

\end{document}